\documentclass[11pt,a4paper]{amsart}
\usepackage[centering,left=3cm,right=3cm]{geometry}
\usepackage{amsfonts,amscd,amssymb,amsmath,amsthm,mathrsfs,xcolor,lscape,longtable}
\usepackage{microtype}
\usepackage[Algorithm]{algorithm}
\usepackage[noend]{algpseudocode}
\usepackage[colorlinks,linkcolor=blue,anchorcolor=blue,citecolor=blue,backref=page]{hyperref}
\usepackage{url}
\usepackage{enumitem}
\usepackage{color}
\usepackage{graphics,epsfig}
\usepackage{float}
\usepackage{longtable}
\usepackage{epstopdf}
\usepackage[utf8]{inputenc}
\usepackage{hyperref}
\usepackage{calc, amscd}
\usepackage[T1]{fontenc}
\usepackage{pdfpages}
\usepackage{mathtools} 
\usepackage{scalerel,stackengine}
\hypersetup{breaklinks=true}
\usepackage{multirow,pbox,lipsum}
\usepackage[thinlines]{easytable}
\usepackage{todonotes}
\usepackage{phaistos}

\usepackage{caption} 

\newtheorem{thm}{Theorem}

\newtheorem{defn}[thm]{Definition}
\newtheorem{cor}[thm]{Corollary}

\theoremstyle{definition}
\newtheorem{rmk}{Remark}

\theoremstyle{definition}
\newtheorem{example}{Example}
\numberwithin{equation}{section}

\newcommand{\Sp}{\mathrm{Sp}}

\renewcommand{\star}{\ast}

\allowdisplaybreaks

\begin{document}
\date{\today}
\title[Thin monodromy in $\mathrm{Sp}(4)$ and $\mathrm{Sp}(6)$]{Thin monodromy in $\mathrm{Sp}(4)$ and $\mathrm{Sp}(6)$}
\author{Jitendra Bajpai, Daniele Dona, Martin Nitsche}
\address{Max-Planck-Institut f\"ur Mathematik, Bonn, Germany}
\email{jitendra@mpim-bonn.mpg.de}
\address{Einstein Institute of Mathematics, Hebrew University of Jerusalem, Israel}
\email{daniele.dona@mail.huji.ac.il}
\address{Institute of Algebra and Geometry, Karlsruhe Institute of Technology, Germany}
\email{martin.nitsche@kit.edu}


\begin{abstract}
We explore the thinness of hypergeometric groups of type $\Sp(4)$ and $\Sp(6)$ by applying a new approach of computer-assisted ping pong.
We prove the thinness of $17$  hypergeometric groups with maximally unipotent monodromy in $\Sp(6)$, completing the classification of all $40$ such groups into arithmetic and thin cases.

In addition, we establish the thinness of further $46$ hypergeometric groups in $\Sp(6)$, and of $3$ hypergeometric groups in $\Sp(4)$, completing the classification of all $\Sp(4)$ hypergeometric groups.
To the best of our knowledge, this article produces the first $63$ examples in the cyclotomic family of Zariski dense non-arithmetic hypergeometric monodromy groups of real rank three.
\end{abstract}

\maketitle
\tableofcontents

\section{Introduction}\label{se:intro}

Let $\theta=z\frac{d}{dz}$ and $\alpha=(\alpha_{1},\ldots,\alpha_{n}),\,\beta=(\beta_{1},\ldots,\beta_{n})\in\mathbb{C}^{n}.$ Then, the \textit{hypergeometric differential equation} of order $n$ is defined by
\begin{equation}\label{hde}
\big[z(\theta+\alpha_{1})\cdot\ldots\cdot(\theta+\alpha_{n})-(\theta+\beta_{1}-1)\cdot\ldots\cdot(\theta+\beta_{n}-1)\big]u(z)=0.
\end{equation}
It is defined on the thrice-punctured Riemann sphere $\mathbb{P}^{1}(\mathbb{C})\setminus\{0,1,\infty\}$, and it has $n$ linearly independent solutions. Thus, the fundamental group $\pi_{1}$ of $\mathbb{P}^{1}(\mathbb{C})\setminus\{0,1,\infty\}$ acts on the (local) solution space of \eqref{hde}: calling $V$ the solution space in the neighbourhood of a point of $\mathbb{P}^{1}(\mathbb{C})\setminus\{0,1,\infty\}$, we get the monodromy representation $\rho:\pi_{1}\rightarrow\mathrm{GL}(V)$. The subgroup $\rho(\pi_{1})$ of $\mathrm{GL}(V)$ is called the \textit{monodromy group} of the hypergeometric differential equation~\eqref{hde}. We also call it the \textit{hypergeometric group} associated to $\alpha,\,\beta\in\mathbb{C}^{n}$.

Levelt (cf.~\cite[Thm.~3.5]{BH}) showed that, if $\alpha_{j}-\beta_{k}\notin\mathbb{Z}$ for all $1\leq j,k\leq n$, then there exists a basis of the solution space of~\eqref{hde} with respect to which the associated hypergeometric group is the subgroup of $\mathrm{GL}_{n}(\mathbb{C})$ generated by the companion matrices $A$ and $B$ of the polynomials
\begin{align*}
f(x) & =\prod_{j=1}^{n}(x-e^{2\pi i\alpha_{j}}), & g(x) & =\prod_{j=1}^{n}(x-e^{2\pi i\beta_{j}})
\end{align*}
respectively, and the monodromy representation is defined by
\begin{align}\label{eq:monodromy}
g_{\infty} & \mapsto A, & g_{0} & \mapsto B^{-1}, & g_{1} & \mapsto A^{-1}B, 
\end{align}
where $g_{0},g_{1},g_{\infty}$ are the loops around $0,1,\infty$ respectively. The condition $\alpha_{j}-\beta_{k}\notin\mathbb{Z}$ for all $1\leq j,k\leq n$ ensures that $f$ and $g$ do not have any common root.

Let $\Gamma(f,g)\subseteq\mathrm{GL}_{n}(\mathbb{C})$ denote the hypergeometric group generated by the companion matrices of the polynomials $f,g$.

\begin{defn}\label{defn:at}
A hypergeometric group $\Gamma(f,g)$ is called arithmetic if it is of finite index in $\mathrm{G}(\mathbb{Z})$, and thin if it has infinite index in $\mathrm{G}(\mathbb{Z})$, where $\mathrm{G}$ is the Zariski closure of $\Gamma(f,g)$ inside $\mathrm{GL}_{n}(\mathbb{C})$.
\end{defn}

Sarnak's question~\cite{Sa14} about classifying the pairs of polynomials $f,g$ for which $\Gamma(f,g)$ is arithmetic or thin has witnessed many interesting developments.
For a detailed discussion on recent progress we refer the interested reader to the introduction of~\cite{bdss}.

In this article we consider the cases $n=4,6$, and we assume that $\alpha,\beta\in\mathbb{Q}^{n}$ and $f,g\in\mathbb{Q}[x]$. This implies that $f,g$ are products of cyclotomic polynomials and that in fact they sit in $\mathbb{Z}[x]$. Up to taking a conjugate of the original group, we may assume $\alpha_{j},\beta_{k}\in[0,1)$, which leads to a list of finitely many cases. Out of these, we only consider the cases that satisfy the following conditions:
\begin{enumerate}
\item $f$ and $g$ have no common roots, such that Levelt's result applies;
\item the Zariski closure of $\Gamma(f,g)$ is $\mathrm{Sp}_\Omega$ for some non-degenerate symplectic form $\Omega$;
\item $\Gamma(f,g)\subseteq\mathrm{Sp}_\Omega$ is primitive (see~\cite[Def.~5.1]{BH}). 
\end{enumerate}

From the results of~\cite{BH} it follows that, up to scalar shift (see~\cite[Def.~5.5]{BH}), there are exactly $58$ such cases in $\Sp(4)$ and $458$ in $\Sp(6)$. For these groups, Sarnak's question becomes whether $\Gamma(f,g)$ has finite or infinite index inside $\mathrm{Sp}_{\Omega}(\mathbb{Z})$.

Of particular relevance are the hypergeometric groups with a \textit{maximally unipotent monodromy}, that is, the hypergeometric groups $\Gamma(f,g)$ with polynomial $f$ associated to $\alpha=\left(0,\ldots,0\right)$. For $n=4$, the 14 hypergeometric groups associated to  $\alpha=\left(0,0,0,0\right)$ emerge as images of monodromy representations arising from Calabi-Yau 3-folds (cf.~\cite[Table~1]{DM06}) and it is expected that many of the 40 groups in the case $n=6$ will similarly arise from Calabi-Yau 5-folds (see for instance~\cite{GMP,LTY}).

The main purpose of this work is to investigate which among the 40 groups with a maximally unipotent monodromy in $\Sp(6)$ are thin. Together with~\cite{bdss, BDN22}, both of which investigate arithmeticity instead, the present article yields a complete classification of these 40 groups. Our methods are effective in tackling other cases in $\Sp(4)$ and $\Sp(6)$ as well.

\subsection{Methods}\label{se:methods}
We adapt the strategy of Brav and Thomas~\cite{BT}: namely, we apply the ping pong lemma to $\langle BA^{-1}\rangle$ and $\langle B\rangle$, where the ping pong table halves are unions of convex cones in $\mathbb{R}^n$. The crucial step, here, is to construct cones for which the procedure works. In dimension $n=4$ Brav and Thomas found a unified formula that applied to all the maximally unipotent cases, from which they produced cones spanned by only $n$ rays, allowing for manual verification.
In dimension $n=6$, in contrast, we use computer assistance to construct the cones. This made it feasible to work with more complicated cones and to search for an individual solution for each case.

To our knowledge, this is the first time that a computer-assisted ping pong proof has been carried out; this is the main methodological novelty of this article. Though not yielding a unified formula, our method has been very successful in practice, leaving out only few cases for which we cannot decide whether they are arithmetic or thin. Furthermore, it seems to us that the method is quite flexible and could also be applied to other groups of matrices. We describe our technique in Section~\ref{se:thin}.

We emphasize that we do not have to rely on high computing power at all. All computations were done on standard consumer hardware, and the crucial verification step, in particular, takes at most a few seconds.

\subsection{Results}

Our main result concerns the $40$ maximally unipotent $\Sp(6)$ hypergeometric groups, listed in Table~\ref{tab:thin} with the same numbering as in~\cite[Table A]{bdss}.

\begin{thm}\label{thm:total}
The degree six symplectic hypergeometric groups with a maximally unipotent monodromy marked as ``thin'' in Table~\ref{tab:thin} are thin. Moreover, they are abstractly isomorphic to the free product $\mathbb{Z}\ast\mathbb{Z}$, except for the cases $\mathrm{A}$-$31$, $\mathrm{A}$-$37$, $\mathrm{A}$-$38$, which are isomorphic to $(\mathbb{Z}\times\mathbb{Z}/2\mathbb{Z})\ast_{\mathbb{Z}/2\mathbb{Z}}\mathbb{Z}/12\mathbb{Z}$, to $\mathbb{Z}\ast\mathbb{Z}/7\mathbb{Z}$, and to $\mathbb{Z}\ast \mathbb{Z}/9\mathbb{Z}$, respectively.
\end{thm}

Together with the arithmeticity results of~\cite{bdss} and~\cite{BDN22}, this gives a complete classification of the maximally unipotent $\Sp(6)$ hypergeometric groups into arithmetic and thin cases.

\begin{cor}
Out of the $40$ degree six symplectic hypergeometric groups with a maximally unipotent monodromy, there are exactly $23$ arithmetic and $17$ thin cases.
\end{cor}

\begin{rmk}
We note that the thin groups of this article are the first examples in the cyclotomic family of Zariski dense non-arithmetic hypergeometric monodromy groups of real rank three.
\end{rmk}

\begin{example}
One interesting example among the 17 thin groups above is the septic case, that is the one with $\beta=\big(\frac{1}{7},\frac{2}{7},\frac{3}{7},\frac{4}{7},\frac{5}{7},\frac{6}{7}\big)$, marked as A-37 in Table~\ref{tab:thin}. This is the second member of the Dwork family to be proved thin, the first being the quintic case of degree four associated to the pair $\alpha=(0,0,0,0)$ and $\beta=\left(\frac{1}{5},\frac{2}{5},\frac{3}{5},\frac{4}{5}\right)$. It is conjectured that all members of the Dwork family are thin: see for instance~\cite[Sect.~3.5]{Sa14}.
\end{example}

\begin{rmk}
After the appearance of the first version of this article on arXiv, Singh and Singh~\cite{SS} were able to independently show the thinness of some of the cases covered in our article. By generalizing to $\mathrm{Sp}(6)$ the unified formula of~\cite{BT}, they obtained solutions that can be checked manually for 7 cases from Table~\ref{tab:thin}, namely A-1 to A-6 and A-10.
\end{rmk}

As a secondary result we address the question of thinness for the remaining $\Sp(6)$ hypergeometric groups, when $\alpha\neq(0,0,0,0,0,0)$.
In Tables~\ref{tab:thin-more} and~\ref{tab:open} below we use the same numbering as in \cite[Table C]{bdss}.

\begin{thm}\label{thm:total-nonmum}
The degree six symplectic hypergeometric groups listed in Table~\ref{tab:thin-more} are thin.
Moreover, these groups are abstractly isomorphic to $\mathbb{Z}\ast\mathbb{Z}$ except for the cases $\mathrm{C}$-$19$, $\mathrm{C}$-$33$, $\mathrm{C}$-$46$, $\mathrm{C}$-$52$, $\mathrm{C}$-$58$ and $\mathrm{C}$-$20$, $\mathrm{C}$-$34$, $\mathrm{C}$-$53$, which are isomorphic to $\mathbb{Z}\ast\mathbb{Z}/7\mathbb{Z}$ and $\mathbb{Z}\ast\mathbb{Z}/9\mathbb{Z}$, respectively.
\end{thm}

Together with the arithmeticity results of~\cite{SV, DFH, bdss, BDN22}, this means that out of the $458$ $\Sp(6)$ hypergeometric groups mentioned in the introduction there are now only $3$ cases whose arithmeticity or thinness is still unknown. We list them in Table~\ref{tab:open}.

Finally, we investigate thinness among hypergeometric groups in $\Sp(4)$. 

\begin{thm}\label{thm:total-sp4}
The degree four symplectic hypergeometric groups labelled as $27$, $35$, $39$ in Table~\ref{tab:sp4-all} are thin.
Moreover, groups $27$ and $35$ are abstractly isomorphic to $\mathbb{Z}\ast\mathbb{Z}$, whereas group $39$ is isomorphic to $\mathbb{Z}\ast\mathbb{Z}/5\mathbb{Z}$.
\end{thm}

Of the $58$ existing cases, $53$ have already been settled in \cite{SV, BT, S15S, S17, BSS}. The result above, together with the arithmeticity result of \cite{BDN22}, completes the classification of all hypergeometric groups in $\Sp(4)$: see Table~\ref{tab:sp4-all}.

\begin{cor}
Out of the $58$ degree four symplectic hypergeometric groups, there are exactly $48$ arithmetic and $10$ thin cases.
\end{cor}

\section{Tables of \texorpdfstring{$\Sp(4)$}{Sp(4)} and \texorpdfstring{$\Sp(6)$}{Sp(6)} hypergeometric groups}
{\centering 
\scriptsize\renewcommand{\arraystretch}{2}
\begin{longtable}{|c|c|c|||c|c|c|}
\caption{Classification of maximally unipotent $\Sp(6)$ hypergeometric groups}\label{tab:thin}\\
\hline
Label & $\beta$ & Nature & Label & $\beta$ & Nature \\
\hline
\hline
A-1 & $\left(\frac{1}{2},\frac{1}{2},\frac{1}{2},\frac{1}{2},\frac{1}{2},\frac{1}{2}\right)$ & thin &
A-2 & $\left(\frac{1}{2},\frac{1}{2},\frac{1}{2},\frac{1}{2},\frac{1}{3},\frac{2}{3}\right)$ & thin \\
\hline
A-3 & $\left(\frac{1}{2},\frac{1}{2},\frac{1}{2},\frac{1}{2},\frac{1}{4},\frac{3}{4}\right)$ & thin &
A-4 & $\left(\frac{1}{2},\frac{1}{2},\frac{1}{2},\frac{1}{2},\frac{1}{6},\frac{5}{6}\right)$ & thin \\
\hline
A-5 & $\left(\frac{1}{2},\frac{1}{2},\frac{1}{3},\frac{1}{3},\frac{2}{3},\frac{2}{3}\right)$ & thin &
A-6 & $\left(\frac{1}{2},\frac{1}{2},\frac{1}{3},\frac{2}{3},\frac{1}{4},\frac{3}{4}\right)$ & thin \\
\hline
A-7 & $\left(\frac{1}{2},\frac{1}{2},\frac{1}{3},\frac{2}{3},\frac{1}{6},\frac{5}{6}\right)$ & thin &
A-8 & $\left(\frac{1}{2},\frac{1}{2},\frac{1}{4},\frac{1}{4},\frac{3}{4},\frac{3}{4}\right)$ & thin \\
\hline
A-9 & $\left(\frac{1}{2},\frac{1}{2},\frac{1}{4},\frac{3}{4},\frac{1}{6},\frac{5}{6}\right)$ & thin &
A-10 & $\left(\frac{1}{2},\frac{1}{2},\frac{1}{5},\frac{2}{5},\frac{3}{5},\frac{4}{5}\right)$ & thin \\
\hline
A-11 & $\left(\frac{1}{2},\frac{1}{2},\frac{1}{6},\frac{1}{6},\frac{5}{6},\frac{5}{6}\right)$ & thin &
A-12 & $\left(\frac{1}{2},\frac{1}{2},\frac{1}{8},\frac{3}{8},\frac{5}{8},\frac{7}{8}\right)$ & thin \\
\hline
A-13 & $\left(\frac{1}{2},\frac{1}{2},\frac{1}{10},\frac{3}{10},\frac{7}{10},\frac{9}{10}\right)$ & thin &
A-14 & $\left(\frac{1}{2},\frac{1}{2},\frac{1}{12},\frac{5}{12},\frac{7}{12},\frac{11}{12}\right)$ & thin \\
\hline
A-15 & $\left(\frac{1}{3},\frac{1}{3},\frac{1}{3},\frac{2}{3},\frac{2}{3},\frac{2}{3}\right)$ & arithmetic \cite{BDN22} &
A-16 & $\left(\frac{1}{3},\frac{1}{3},\frac{2}{3},\frac{2}{3},\frac{1}{4},\frac{3}{4}\right)$ & arithmetic \cite{BDN22} \\
\hline
A-17 & $\left(\frac{1}{3},\frac{1}{3},\frac{2}{3},\frac{2}{3},\frac{1}{6},\frac{5}{6}\right)$ & arithmetic \cite{bdss} &
A-18 & $\left(\frac{1}{3},\frac{2}{3},\frac{1}{4},\frac{1}{4},\frac{3}{4},\frac{3}{4}\right)$ & arithmetic \cite{bdss} \\
\hline
A-19 & $\left(\frac{1}{3},\frac{2}{3},\frac{1}{4},\frac{3}{4},\frac{1}{6},\frac{5}{6}\right)$ & arithmetic \cite{bdss} &
A-20 & $\left(\frac{1}{3},\frac{2}{3},\frac{1}{6},\frac{5}{6},\frac{1}{6},\frac{5}{6}\right)$ & arithmetic \cite{bdss} \\
\hline
A-21 & $\left(\frac{1}{3},\frac{2}{3},\frac{1}{5},\frac{2}{5},\frac{3}{5},\frac{4}{5}\right)$ & arithmetic \cite{BDN22} &
A-22 & $\left(\frac{1}{3},\frac{2}{3},\frac{1}{8},\frac{3}{8},\frac{5}{8},\frac{7}{8}\right)$ & arithmetic \cite{bdss} \\
\hline
A-23 & $\left(\frac{1}{3},\frac{2}{3},\frac{1}{10},\frac{3}{10},\frac{7}{10},\frac{9}{10}\right)$ & arithmetic \cite{bdss} &
A-24 & $\left(\frac{1}{3},\frac{2}{3},\frac{1}{12},\frac{5}{12},\frac{7}{12},\frac{11}{12}\right)$ & arithmetic \cite{BDN22} \\
\hline
A-25 & $\left(\frac{1}{4},\frac{1}{4},\frac{1}{4},\frac{3}{4},\frac{3}{4},\frac{3}{4}\right)$ & arithmetic \cite{bdss} &
A-26 & $\left(\frac{1}{4},\frac{1}{4},\frac{3}{4},\frac{3}{4},\frac{1}{6},\frac{5}{6}\right)$ & arithmetic \cite{bdss} \\
\hline
A-27 & $\left(\frac{1}{4},\frac{3}{4},\frac{1}{5},\frac{2}{5},\frac{3}{5},\frac{4}{5}\right)$ & arithmetic \cite{bdss} &
A-28 & $\left(\frac{1}{4},\frac{3}{4},\frac{1}{6},\frac{5}{6},\frac{1}{6},\frac{5}{6}\right)$ & arithmetic \cite{bdss} \\
\hline
A-29 & $\left(\frac{1}{4},\frac{3}{4},\frac{1}{8},\frac{3}{8},\frac{5}{8},\frac{7}{8}\right)$ & arithmetic \cite{bdss} &
A-30 & $\left(\frac{1}{4},\frac{3}{4},\frac{1}{10},\frac{3}{10},\frac{7}{10},\frac{9}{10}\right)$ & arithmetic \cite{bdss} \\
\hline
A-31 & $\left(\frac{1}{4},\frac{3}{4},\frac{1}{12},\frac{5}{12},\frac{7}{12},\frac{11}{12}\right)$ & thin &
A-32 & $\left(\frac{1}{5},\frac{2}{5},\frac{3}{5},\frac{4}{5},\frac{1}{6},\frac{5}{6}\right)$ & arithmetic \cite{bdss} \\
\hline
A-33 & $\left(\frac{1}{6},\frac{5}{6},\frac{1}{6},\frac{5}{6},\frac{1}{6},\frac{5}{6}\right)$ & arithmetic \cite{bdss} &
A-34 & $\left(\frac{1}{6},\frac{5}{6},\frac{1}{8},\frac{3}{8},\frac{5}{8},\frac{7}{8}\right)$ & arithmetic \cite{bdss} \\
\hline
A-35 & $\left(\frac{1}{6},\frac{5}{6},\frac{1}{10},\frac{3}{10},\frac{7}{10},\frac{9}{10}\right)$ & arithmetic \cite{bdss} &
A-36 & $\left(\frac{1}{6},\frac{5}{6},\frac{1}{12},\frac{5}{12},\frac{7}{12},\frac{11}{12}\right)$ & arithmetic \cite{bdss} \\
\hline
A-37 & $\left(\frac{1}{7},\frac{2}{7},\frac{3}{7},\frac{4}{7},\frac{5}{7},\frac{6}{7}\right)$ & thin &
A-38 & $\left(\frac{1}{9},\frac{2}{9},\frac{4}{9},\frac{5}{9},\frac{7}{9},\frac{8}{9}\right)$ & thin \\
\hline
A-39 & $\left(\frac{1}{14},\frac{3}{14},\frac{5}{14},\frac{9}{14},\frac{11}{14},\frac{13}{14}\right)$ & arithmetic \cite{BDN22} &
A-40 & $\left(\frac{1}{18},\frac{5}{18},\frac{7}{18},\frac{11}{18},\frac{13}{18},\frac{17}{18}\right)$ & arithmetic \cite{bdss} \\
\hline
\end{longtable}}

{\centering
\scriptsize\renewcommand{\arraystretch}{2}
\begin{longtable}{|c|c|c|||c|c|c|}
\caption{More thinness by ping-pong: $46$ new examples of $\Sp(6)$ hypergeometric groups whose thinness follows from Sections~\ref{se:thin} and~\ref{se:morethin}.}\label{tab:thin-more}\\
\hline
 Label & $\alpha$ & $\beta$ & Label & $\alpha$ & $\beta$   \\                                                            
\hline
\hline
C-2 & $(0,0,0,0,\frac{1}{3},\frac{2}{3})$ &  $\big(\frac{1}{2},\frac{1}{2},\frac{1}{2},\frac{1}{2},\frac{1}{4},\frac{3}{4}\big)$ &  
C-3 & $(0,0,0,0,\frac{1}{3},\frac{2}{3})$ & $\big(\frac{1}{2},\frac{1}{2},\frac{1}{2},\frac{1}{2},\frac{1}{6},\frac{5}{6}\big)$  \\
\hline
C-4 & $\left(0,0,0,0,\frac{1}{3},\frac{2}{3}\right)$&$\left(\frac{1}{2},\frac{1}{2},\frac{1}{4},\frac{1}{4},\frac{3}{4},\frac{3}{4}\right)$ & 
C-5 & $(0,0,0,0,\frac{1}{3},\frac{2}{3})$ & $\big(\frac{1}{2},\frac{1}{2},\frac{1}{5},\frac{2}{5},\frac{3}{5},\frac{4}{5}\big)$ \\
\hline
C-6 & $(0,0,0,0,\frac{1}{3},\frac{2}{3})$ & $\big(\frac{1}{2},\frac{1}{2},\frac{1}{8},\frac{3}{8},\frac{5}{8},\frac{7}{8}\big)$ & 
C-7 & $(0,0,0,0,\frac{1}{3},\frac{2}{3})$&$\big(\frac{1}{2},\frac{1}{2},\frac{1}{10},\frac{3}{10},\frac{7}{10},\frac{9}{10}\big)$ \\
\hline
C-8 & $(0,0,0,0,\frac{1}{3},\frac{2}{3})$ & $\big(\frac{1}{2},\frac{1}{2},\frac{1}{12},\frac{5}{12},\frac{7}{12},\frac{11}{12}\big)$ & 
C-11 & $(0,0,0,0,\frac{1}{4},\frac{3}{4})$ & $\big(\frac{1}{2},\frac{1}{2},\frac{1}{2},\frac{1}{2},\frac{1}{3},\frac{2}{3}\big)$ \\
\hline
C-12 & $\left(0,0,0,0,\frac{1}{4},\frac{3}{4}\right)$ & $\left(\frac{1}{2},\frac{1}{2},\frac{1}{3},\frac{1}{3},\frac{2}{3},\frac{2}{3}\right)$ & 
C-13 & $(0,0,0,0,\frac{1}{4},\frac{3}{4})$ & $\big(\frac{1}{2},\frac{1}{2},\frac{1}{3},\frac{2}{3},\frac{1}{6},\frac{5}{6}\big)$  \\
\hline
C-14 & $(0,0,0,0,\frac{1}{4},\frac{3}{4})$ & $\big(\frac{1}{2},\frac{1}{2},\frac{1}{5},\frac{2}{5},\frac{3}{5},\frac{4}{5}\big)$ & 
C-15 & $\left(0,0,0,0,\frac{1}{4},\frac{3}{4}\right)$ & $\left(\frac{1}{2},\frac{1}{2},\frac{1}{6},\frac{5}{6},\frac{5}{6},\frac{5}{6}\right)$ \\
\hline
C-16 & $(0,0,0,0,\frac{1}{4},\frac{3}{4})$ & $\big(\frac{1}{2},\frac{1}{2},\frac{1}{8},\frac{3}{8},\frac{5}{8},\frac{7}{8}\big)$  & 
C-17 & $(0,0,0,0,\frac{1}{4},\frac{3}{4})$& $\big(\frac{1}{2},\frac{1}{2},\frac{1}{10},\frac{3}{10},\frac{7}{10},\frac{9}{10}\big)$ \\
\hline
C-18 &  $(0,0,0,0,\frac{1}{4},\frac{3}{4})$& $\big(\frac{1}{2},\frac{1}{2},\frac{1}{12},\frac{5}{12},\frac{7}{12},\frac{11}{12}\big)$ & 
C-19 & $\left(0,0,0,0,\frac{1}{4},\frac{3}{4}\right)$&$\left(\frac{1}{7},\frac{2}{7},\frac{3}{7},\frac{4}{7},\frac{5}{7},\frac{6}{7}\right)$\\
\hline
C-20 & $(0,0,0,0,\frac{1}{4},\frac{3}{4})$& $\big(\frac{1}{9},\frac{2}{9},\frac{4}{9},\frac{5}{9},\frac{7}{9},\frac{8}{9}\big)$ & 
C-21 & $(0,0,0,0,\frac{1}{6},\frac{5}{6})$  & $\big(\frac{1}{2},\frac{1}{2},\frac{1}{2},\frac{1}{2},\frac{1}{3},\frac{2}{3}\big)$\\
\hline
C-22 & $\left(0,0,0,0,\frac{1}{6},\frac{5}{6}\right)$ & $\left(\frac{1}{2},\frac{1}{2},\frac{1}{3},\frac{1}{3},\frac{2}{3},\frac{2}{3}\right)$ & 
C-23 & $(0,0,0,0,\frac{1}{6},\frac{5}{6})$ & $\big(\frac{1}{2},\frac{1}{2},\frac{1}{3},\frac{2}{3},\frac{1}{4},\frac{3}{4}\big)$  \\
\hline
C-24 & $\left(0,0,0,0,\frac{1}{6},\frac{5}{6}\right)$ & $\left(\frac{1}{2},\frac{1}{2},\frac{1}{4},\frac{1}{4},\frac{3}{4},\frac{3}{4}\right)$ & 
C-25 & $(0,0,0,0,\frac{1}{6},\frac{5}{6})$&  $\big(\frac{1}{2},\frac{1}{2},\frac{1}{5},\frac{2}{5},\frac{3}{5},\frac{4}{5}\big)$ \\
\hline
C-26 & $(0,0,0,0,\frac{1}{6},\frac{5}{6})$& $\big(\frac{1}{2},\frac{1}{2},\frac{1}{8},\frac{3}{8},\frac{5}{8},\frac{7}{8}\big)$   & 
C-27 & $(0,0,0,0,\frac{1}{6},\frac{5}{6})$& $\big(\frac{1}{2},\frac{1}{2},\frac{1}{10},\frac{3}{10},\frac{7}{10},\frac{9}{10}\big)$  \\
\hline
C-28 & $(0,0,0,0,\frac{1}{6},\frac{5}{6})$ & $\big(\frac{1}{2},\frac{1}{2},\frac{1}{12},\frac{5}{12},\frac{7}{12},\frac{11}{12}\big)$  & 
C-33 & $(0,0,0,0,\frac{1}{6},\frac{5}{6})$ &  $\big(\frac{1}{7},\frac{2}{7},\frac{3}{7},\frac{4}{7},\frac{5}{7},\frac{6}{7}\big)$ \\ 
\hline 
C-34 & $(0,0,0,0,\frac{1}{6},\frac{5}{6})$& $\big(\frac{1}{9},\frac{2}{9},\frac{4}{9},\frac{5}{9},\frac{7}{9},\frac{8}{9}\big)$   & 
C-35 & $(0,0,\frac{1}{3},\frac{2}{3},\frac{1}{4},\frac{3}{4})$  & $(\frac{1}{2},\frac{1}{2},\frac{1}{5},\frac{2}{5},\frac{3}{5},\frac{4}{5})$ \\
\hline
C-36 & $(0,0,\frac{1}{3},\frac{2}{3},\frac{1}{6},\frac{5}{6})$ & $(\frac{1}{2},\frac{1}{2},\frac{1}{4},\frac{1}{4},\frac{3}{4},\frac{3}{4})$  & 
C-37 & $(0,0,\frac{1}{3},\frac{2}{3},\frac{1}{6},\frac{5}{6})$  & $(\frac{1}{2},\frac{1}{2},\frac{1}{5},\frac{2}{5},\frac{3}{5},\frac{4}{5})$ \\
\hline
C-38 & $(0,0,\frac{1}{3},\frac{2}{3},\frac{1}{6},\frac{5}{6})$ & $(\frac{1}{2},\frac{1}{2},\frac{1}{8},\frac{3}{8},\frac{5}{8},\frac{7}{8})$ & 
C-40 & $\left(0,0,\frac{1}{4},\frac{1}{4},\frac{3}{4},\frac{3}{4}\right)$&$\left(\frac{1}{2},\frac{1}{2},\frac{1}{3},\frac{2}{3},\frac{1}{3},\frac{2}{3}\right)$\\
\hline
C-41 & $\left(0,0,\frac{1}{4},\frac{1}{4},\frac{3}{4},\frac{3}{4}\right)$& $\left(\frac{1}{2},\frac{1}{2},\frac{1}{5},\frac{2}{5},\frac{3}{5},\frac{4}{5}\right)$ & 
C-43 &$\left(0,0,\frac{1}{4},\frac{3}{4},\frac{1}{6},\frac{5}{6}\right)$ & $\left(\frac{1}{2},\frac{1}{2},\frac{1}{3},\frac{1}{3},\frac{2}{3},\frac{2}{3}\right)$ \\
\hline
C-44 & $(0,0,\frac{1}{4},\frac{3}{4},\frac{1}{6},\frac{5}{6})$&  $(\frac{1}{2},\frac{1}{2},\frac{1}{5},\frac{2}{5},\frac{3}{5},\frac{4}{5})$ & 
C-45 & $(0,0,\frac{1}{4},\frac{3}{4},\frac{1}{6},\frac{5}{6})$& $(\frac{1}{2},\frac{1}{2},\frac{1}{8},\frac{3}{8},\frac{5}{8},\frac{7}{8})$  \\
\hline
C-46 & $\left(0,0,\frac{1}{4},\frac{3}{4},\frac{1}{6},\frac{5}{6}\right)$ & $\left(\frac{1}{7},\frac{2}{7},\frac{3}{7},\frac{4}{7},\frac{5}{7},\frac{6}{7}\right)$ &  
C-48 & $\left(0,0,\frac{1}{6},\frac{1}{6},\frac{5}{6},\frac{5}{6}\right)$ & $\left(\frac{1}{2},\frac{1}{2},\frac{1}{3},\frac{1}{3},\frac{2}{3},\frac{2}{3}\right)$ \\
\hline
C-49 & $(0,0,\frac{1}{6},\frac{1}{6},\frac{5}{6},\frac{5}{6})$ &  $(\frac{1}{2},\frac{1}{2},\frac{1}{5},\frac{2}{5},\frac{3}{5},\frac{4}{5})$ & 
C-50 & $(0,0,\frac{1}{6},\frac{1}{6},\frac{5}{6},\frac{5}{6})$ & $(\frac{1}{2},\frac{1}{2},\frac{1}{8},\frac{3}{8},\frac{5}{8},\frac{7}{8})$  \\
\hline
C-52 & $(0,0,\frac{1}{6},\frac{1}{6},\frac{5}{6},\frac{5}{6})$&$(\frac{1}{7},\frac{2}{7},\frac{3}{7},\frac{4}{7},\frac{5}{7},\frac{6}{7})$  & 
C-53  & $\left(0,0,\frac{1}{6},\frac{1}{6},\frac{5}{6},\frac{5}{6}\right)$ & $\left(\frac{1}{9},\frac{2}{9},\frac{4}{9},\frac{5}{9},\frac{7}{9},\frac{8}{9}\right)$\\
\hline
C-54 & $(0,0,\frac{1}{8},\frac{3}{8},\frac{5}{8},\frac{7}{8})$&$(\frac{1}{2},\frac{1}{2},\frac{1}{5},\frac{2}{5},\frac{3}{5},\frac{4}{5})$ & 
C-56 & $(0,0,\frac{1}{10},\frac{3}{10},\frac{7}{10},\frac{9}{10})$  & $(\frac{1}{2},\frac{1}{2},\frac{1}{5},\frac{2}{5},\frac{3}{5},\frac{4}{5})$ \\
\hline
C-57 & $(0,0,\frac{1}{10},\frac{3}{10},\frac{7}{10},\frac{9}{10})$ & $(\frac{1}{2},\frac{1}{2},\frac{1}{12},\frac{5}{12},\frac{7}{12},\frac{11}{12})$ &
C-58 & $(0,0,\frac{1}{10},\frac{3}{10},\frac{7}{10},\frac{9}{10})$ & $(\frac{1}{7},\frac{2}{7},\frac{3}{7},\frac{4}{7},\frac{5}{7},\frac{6}{7})$ \\
\hline
\end{longtable}}

{\centering
\scriptsize\renewcommand{\arraystretch}{2}
\begin{longtable}{|c|c|c|||c|c|c|}
\caption{Three cases of $\Sp(6)$ hypergeometric groups whose arithmeticity or thinness is still \emph{unknown}.}\label{tab:open}\\\hline
Label & $\alpha$ & $\beta$ & Label & $\alpha$ & $\beta$\\\hline\hline
C-32 & $\left(0,0,0,0,\frac{1}{6},\frac{5}{6}\right)$&$\left(\frac{1}{4},\frac{3}{4},\frac{1}{12},\frac{5}{12},\frac{7}{12},\frac{11}{12}\right)$ &
C-47 &$\left(0,0,\frac{1}{5},\frac{2}{5},\frac{3}{5},\frac{4}{5}\right)$ & $\left(\frac{1}{2},\frac{1}{2},\frac{1}{3},\frac{1}{3},\frac{2}{3},\frac{2}{3}\right)$ \\\hline
C-55 &$\left(0,0,\frac{1}{8},\frac{3}{8},\frac{5}{8},\frac{7}{8}\right)$ & $\left(\frac{1}{2},\frac{1}{2},\frac{1}{12},\frac{5}{12},\frac{7}{12},\frac{11}{12}\right)$ & & &\\\hline
\end{longtable}}

\newcounter{rownum-0}
\setcounter{rownum-0}{0}
{\centering
\scriptsize\renewcommand{\arraystretch}{2}
\begin{longtable}{|c|c|c|c||c|c|c|c|c|c|}
\caption{Classification of $\Sp(4)$ hypergeometric groups}\label{tab:sp4-all}\\\hline
Label & $\alpha$  & $\beta$  & Nature & Label & $\alpha$  & $\beta$  & Nature   \\
\hline
\hline
\refstepcounter{rownum-0}\arabic{rownum-0} &  $(0,0,0,0)$  & $\big(\frac{1}{2},\frac{1}{2},\frac{1}{2},\frac{1}{2}\big)$ & thin~\cite{BT}& \refstepcounter{rownum-0}\arabic{rownum-0} &  $(0,0,0,0)$ & $\big(\frac{1}{2},\frac{1}{2},\frac{1}{3},\frac{2}{3}\big)$& thin~\cite{BT} \\
\hline
\refstepcounter{rownum-0}\arabic{rownum-0} &  $(0,0,0,0)$  & $\big(\frac{1}{2},\frac{1}{2},\frac{1}{4},\frac{3}{4}\big)$ & thin~\cite{BT} & \refstepcounter{rownum-0}\arabic{rownum-0} &  $(0,0,0,0)$ & $\big(\frac{1}{2},\frac{1}{2},\frac{1}{6},\frac{5}{6}\big)$ & thin~\cite{BT}\\
\hline
\refstepcounter{rownum-0}\arabic{rownum-0} &  $(0,0,0,0)$  & $\big(\frac{1}{3},\frac{1}{3},\frac{2}{3},\frac{2}{3}\big)$& arithmetic~\cite{S15S}  & \refstepcounter{rownum-0}\arabic{rownum-0} &  $(0,0,0,0)$ & $\big(\frac{1}{3},\frac{2}{3},\frac{1}{4},\frac{3}{4}\big)$  & arithmetic~\cite{S15S}   \\
\hline
\refstepcounter{rownum-0}\arabic{rownum-0} &  $(0,0,0,0)$  & $\big(\frac{1}{3},\frac{2}{3},\frac{1}{6},\frac{5}{6}\big)$ & arithmetic~\cite{S15S}   & \refstepcounter{rownum-0}\arabic{rownum-0} &  $(0,0,0,0)$ & $\big(\frac{1}{4},\frac{1}{4},\frac{3}{4},\frac{3}{4}\big)$ & arithmetic~\cite{S15S}   \\
\hline
\refstepcounter{rownum-0}\arabic{rownum-0} &  $(0,0,0,0)$  & $\big(\frac{1}{4},\frac{3}{4},\frac{1}{6},\frac{5}{6}\big)$ & arithmetic~\cite{SV}  & \refstepcounter{rownum-0}\arabic{rownum-0} &  $(0,0,0,0)$ & $\big(\frac{1}{5},\frac{2}{5},\frac{3}{5},\frac{4}{5}\big)$ & thin~\cite{BT}\\
\hline
\refstepcounter{rownum-0}\arabic{rownum-0} &  $(0,0,0,0)$  &  $\big(\frac{1}{6},\frac{1}{6},\frac{5}{6},\frac{5}{6}\big)$ & arithmetic~\cite{SV} & \refstepcounter{rownum-0}\arabic{rownum-0} &  $(0,0,0,0)$ & $\big(\frac{1}{8},\frac{3}{8},\frac{5}{8},\frac{7}{8}\big)$ & thin~\cite{BT}\\
\hline
\refstepcounter{rownum-0}\arabic{rownum-0} &  $(0,0,0,0)$  & $\big(\frac{1}{10},\frac{3}{10},\frac{7}{10},\frac{9}{10}\big)$ & arithmetic~\cite{SV} & \refstepcounter{rownum-0}\arabic{rownum-0} &  $(0,0,0,0)$ & $\big(\frac{1}{12},\frac{5}{12},\frac{7}{12},\frac{11}{12}\big)$ & thin~\cite{BT} \\
\hline
\refstepcounter{rownum-0}\arabic{rownum-0}\label{A8} & $(0,0,\frac{1}{2},\frac{1}{2})$ & $\big(\frac{1}{3},\frac{1}{3},\frac{2}{3},\frac{2}{3}\big)$ & arithmetic~\cite{SV} & \refstepcounter{rownum-0}\arabic{rownum-0} &  $(0,0,\frac{1}{2},\frac{1}{2})$ & $\big(\frac{1}{3},\frac{2}{3},\frac{1}{4},\frac{3}{4}\big)$ & arithmetic~\cite{SV}\\
\hline
\refstepcounter{rownum-0}\arabic{rownum-0}\label{A} & $(0,0,\frac{1}{2},\frac{1}{2})$ & $ \big(\frac{1}{5},\frac{2}{5},\frac{3}{5},\frac{4}{5}\big)$ & arithmetic~\cite{SV} & \refstepcounter{rownum-0}\arabic{rownum-0} &  $(0,0,\frac{1}{3},\frac{2}{3})$ & $\big(\frac{1}{2},\frac{1}{2},\frac{1}{4},\frac{3}{4}\big)$ & arithmetic~\cite{BSS} \\
\hline
\refstepcounter{rownum-0}\arabic{rownum-0} & $(0,0,\frac{1}{3},\frac{2}{3})$ & $\big(\frac{1}{2},\frac{1}{2},\frac{1}{6},\frac{5}{6}\big)$ & arithmetic~\cite{SV} & \refstepcounter{rownum-0}\arabic{rownum-0} &  $(0,0,\frac{1}{3},\frac{2}{3})$ &  $\big(\frac{1}{4},\frac{1}{4},\frac{3}{4},\frac{3}{4}\big)$ & arithmetic~\cite{SV} \\
\hline
\refstepcounter{rownum-0}\arabic{rownum-0} & $(0,0,\frac{1}{3},\frac{2}{3})$ & $\big(\frac{1}{4},\frac{3}{4},\frac{1}{6},\frac{5}{6}\big)$  & arithmetic~\cite{SV} & \refstepcounter{rownum-0}\arabic{rownum-0} &  $(0,0,\frac{1}{3},\frac{2}{3})$ & $\big(\frac{1}{5},\frac{2}{5},\frac{3}{5},\frac{4}{5} \big)$ & arithmetic~\cite{SV} \\
\hline
\refstepcounter{rownum-0}\arabic{rownum-0} & $(0,0,\frac{1}{3},\frac{2}{3})$ & $\big(\frac{1}{6},\frac{5}{6},\frac{1}{6},\frac{5}{6}\big)$ & arithmetic~\cite{SV} & \refstepcounter{rownum-0}\arabic{rownum-0} &  $(0,0,\frac{1}{3},\frac{2}{3})$ & $\big(\frac{1}{8},\frac{3}{8},\frac{5}{8},\frac{7}{8}\big)$ & arithmetic~\cite{SV} \\
\hline
\refstepcounter{rownum-0}\arabic{rownum-0} & $(0,0,\frac{1}{3},\frac{2}{3})$ & $\big(\frac{1}{10},\frac{3}{10},\frac{7}{10},\frac{9}{10} \big)$ & arithmetic~\cite{SV} & \refstepcounter{rownum-0}\arabic{rownum-0} &  $(0,0,\frac{1}{3},\frac{2}{3})$ & $\big( \frac{1}{12},\frac{5}{12},\frac{7}{12},\frac{11}{12}\big)$ & arithmetic~\cite{SV} \\
\hline
\refstepcounter{rownum-0}\arabic{rownum-0} & $(0,0,\frac{1}{4},\frac{3}{4})$ & $ \big(\frac{1}{2},\frac{1}{2},\frac{1}{3},\frac{2}{3}\big)$   & thin & \refstepcounter{rownum-0}\arabic{rownum-0} &   $(0,0,\frac{1}{4},\frac{3}{4})$ & $\big(\frac{1}{3},\frac{1}{3},\frac{2}{3},\frac{2}{3}\big)$  & arithmetic~\cite{S17}\\
\hline
\refstepcounter{rownum-0}\arabic{rownum-0} & $(0,0,\frac{1}{4},\frac{3}{4})$ & $\big( \frac{1}{3},\frac{2}{3},\frac{1}{6},\frac{5}{6}\big)$ & arithmetic~\cite{SV} & \refstepcounter{rownum-0}\arabic{rownum-0} &  $(0,0,\frac{1}{4},\frac{3}{4})$ & $\left(\frac{1}{5},\frac{2}{5},\frac{3}{5},\frac{4}{5}\right)$ & arithmetic~\cite{BDN22} \\
\hline
\refstepcounter{rownum-0}\arabic{rownum-0} & $(0,0,\frac{1}{4},\frac{3}{4})$ &$\big( \frac{1}{6},\frac{5}{6},\frac{1}{6},\frac{5}{6}\big)$ & arithmetic~\cite{SV} & \refstepcounter{rownum-0}\arabic{rownum-0} &  $(0,0,\frac{1}{4},\frac{3}{4})$ &  $\big(\frac{1}{8},\frac{3}{8},\frac{5}{8},\frac{7}{8} \big)$& arithmetic~\cite{SV}\\
\hline
\refstepcounter{rownum-0}\arabic{rownum-0} & $(0,0,\frac{1}{4},\frac{3}{4})$ &$\big(\frac{1}{10},\frac{3}{10},\frac{7}{10},\frac{9}{10} \big)$ & arithmetic~\cite{SV} & \refstepcounter{rownum-0}\arabic{rownum-0} &  $(0,0,\frac{1}{4},\frac{3}{4})$ &  $\big(\frac{1}{12},\frac{5}{12},\frac{7}{12},\frac{11}{12}\big)$& arithmetic~\cite{SV} \\
\hline
\refstepcounter{rownum-0}\arabic{rownum-0}\label{T9} & $(0,0,\frac{1}{6},\frac{5}{6})$ &$\big(\frac{1}{2},\frac{1}{2},\frac{1}{3},\frac{2}{3}\big)$& thin & 
\refstepcounter{rownum-0}\arabic{rownum-0} & $(0,0,\frac{1}{6},\frac{5}{6})$ & $\big(\frac{1}{3},\frac{1}{3},\frac{2}{3},\frac{2}{3}\big)$ &  arithmetic~\cite{S17} \\
\hline
\refstepcounter{rownum-0}\arabic{rownum-0} & $(0,0,\frac{1}{6},\frac{5}{6})$ &  $\big(\frac{1}{3},\frac{2}{3},\frac{1}{4},\frac{3}{4}\big)$ & arithmetic~\cite{S17} & 
\refstepcounter{rownum-0}\arabic{rownum-0} & $(0,0,\frac{1}{6},\frac{5}{6})$ & $\big( \frac{1}{4},\frac{1}{4},\frac{3}{4},\frac{3}{4}\big)$ & arithmetic~\cite{S17}  \\
\hline
\refstepcounter{rownum-0}\arabic{rownum-0} & $(0,0,\frac{1}{6},\frac{5}{6})$ &$\big(\frac{1}{5},\frac{2}{5},\frac{3}{5},\frac{4}{5} \big)$ & thin & 
\refstepcounter{rownum-0}\arabic{rownum-0} & $(0,0,\frac{1}{6},\frac{5}{6})$ & $\big( \frac{1}{8},\frac{3}{8},\frac{5}{8},\frac{7}{8}\big)$  &  arithmetic~\cite{BDN22} \\
\hline
\refstepcounter{rownum-0}\arabic{rownum-0} & $(0,0,\frac{1}{6},\frac{5}{6})$ & $\big(\frac{1}{10},\frac{3}{10},\frac{7}{10},\frac{9}{10}  \big)$ & arithmetic~\cite{SV} & 
\refstepcounter{rownum-0}\arabic{rownum-0} & $(0,0,\frac{1}{6},\frac{5}{6})$ & $\big( \frac{1}{12},\frac{5}{12},\frac{7}{12},\frac{11}{12} \big)$  &  arithmetic~\cite{S17} \\
\hline
\refstepcounter{rownum-0}\arabic{rownum-0} &$(\frac{1}{3},\frac{1}{3}, \frac{2}{3},\frac{2}{3})$ &  $\big(  \frac{1}{4},\frac{1}{4},\frac{3}{4},\frac{3}{4}\big)$   & arithmetic~\cite{SV} &
 \refstepcounter{rownum-0}\arabic{rownum-0} & $( \frac{1}{3},\frac{1}{3}, \frac{2}{3},\frac{2}{3})$ & $\big( \frac{1}{4},\frac{3}{4},\frac{1}{6},\frac{5}{6} \big)$ & arithmetic~\cite{S17}\\
\hline
 \refstepcounter{rownum-0}\arabic{rownum-0} & $(\frac{1}{3},\frac{1}{3}, \frac{2}{3},\frac{2}{3})$ & $\big(  \frac{1}{5},\frac{2}{5},\frac{3}{5},\frac{4}{5} \big)$  & arithmetic~\cite{SV}&
 \refstepcounter{rownum-0}\arabic{rownum-0} &$(\frac{1}{3},\frac{1}{3}, \frac{2}{3},\frac{2}{3})$ & $\big( \frac{1}{6},\frac{1}{6} ,\frac{5}{6},\frac{5}{6}  \big)$& arithmetic~\cite{S17} \\
\hline  
 \refstepcounter{rownum-0}\arabic{rownum-0} &$(\frac{1}{3},\frac{1}{3}, \frac{2}{3},\frac{2}{3})$ & $\big(  \frac{1}{8},\frac{3}{8},\frac{5}{8},\frac{7}{8}\big)$ & arithmetic~\cite{SV}&
 \refstepcounter{rownum-0}\arabic{rownum-0} & $(\frac{1}{3},\frac{1}{3}, \frac{2}{3},\frac{2}{3})$ & $\big( \frac{1}{10},\frac{3}{10},\frac{7}{10},\frac{9}{10}  \big)$ & arithmetic~\cite{S17}\\
\hline
 \refstepcounter{rownum-0}\arabic{rownum-0} & $(\frac{1}{3},\frac{1}{3}, \frac{2}{3},\frac{2}{3})$ & $\big( \frac{1}{12},\frac{5}{12},\frac{7}{12},\frac{11}{12}  \big)$ &  arithmetic~\cite{SV} &
 \refstepcounter{rownum-0}\arabic{rownum-0} & $(\frac{1}{3}, \frac{2}{3},\frac{1}{4},\frac{3}{4})$ & $\big(  \frac{1}{5},\frac{2}{5},\frac{3}{5},\frac{4}{5} \big)$ & arithmetic~\cite{SV}\\
\hline
 \refstepcounter{rownum-0}\arabic{rownum-0} & $(\frac{1}{3}, \frac{2}{3},\frac{1}{4},\frac{3}{4})$ & $\big( \frac{1}{8},\frac{3}{8},\frac{5}{8},\frac{7}{8} \big)$& arithmetic~\cite{SV} &
 \refstepcounter{rownum-0}\arabic{rownum-0} &$(\frac{1}{3}, \frac{2}{3},\frac{1}{4},\frac{3}{4})$ & $\big( \frac{1}{10},\frac{3}{10},\frac{7}{10},\frac{9}{10}  \big)$ & arithmetic~\cite{SV}\\
\hline 
 \refstepcounter{rownum-0}\arabic{rownum-0} & $(\frac{1}{3}, \frac{2}{3},\frac{1}{4},\frac{3}{4})$ & $\big( \frac{1}{12},\frac{5}{12},\frac{7}{12},\frac{11}{12}  \big)$ & arithmetic~\cite{SV}&
 \refstepcounter{rownum-0}\arabic{rownum-0} & $(\frac{1}{3}, \frac{2}{3},\frac{1}{6},\frac{5}{6})$ & $\big(  \frac{1}{5},\frac{2}{5},\frac{3}{5},\frac{4}{5}  \big)$ &  arithmetic~\cite{SV} \\
\hline 
 \refstepcounter{rownum-0}\arabic{rownum-0}&$(\frac{1}{4},\frac{1}{4},\frac{3}{4},\frac{3}{4})$ & $\big(\frac{1}{5},\frac{2}{5},\frac{3}{5},\frac{4}{5}   \big)$ & arithmetic~\cite{SV} &
   \refstepcounter{rownum-0}\arabic{rownum-0} & $( \frac{1}{5},\frac{2}{5},\frac{3}{5},\frac{4}{5} )$ & $\big( \frac{1}{8},\frac{3}{8},\frac{5}{8},\frac{7}{8} \big)$ & arithmetic~\cite{SV} \\
\hline 
 \refstepcounter{rownum-0}\arabic{rownum-0}& $( \frac{1}{5},\frac{2}{5},\frac{3}{5},\frac{4}{5})$ & $\big( \frac{1}{10},\frac{3}{10},\frac{7}{10},\frac{9}{10} \big)$ & arithmetic~\cite{SV} &
 \refstepcounter{rownum-0}\arabic{rownum-0} & $(\frac{1}{5},\frac{2}{5},\frac{3}{5},\frac{4}{5})$ &$\big( \frac{1}{12},\frac{5}{12},\frac{7}{12},\frac{11}{12}  \big)$  & arithmetic~\cite{SV}\\
\hline
\end{longtable}}


\section{Ping pong: methodology and code}\label{se:thin}

\subsection{Ping-pong}\label{se:pingpong}

In the present section we discuss in more detail the ping pong method. Our discussion focuses on $\mathrm{Sp}(6)$ since our main result is Theorem~\ref{thm:total}, but the method can be effectively used in more generality.

To show that the groups in Tables~\ref{tab:thin} and~\ref{tab:thin-more} are thin, we will actually prove a stronger result. Namely, let $A$ and $B$ be the companion matrices generating the group $G$, as in the introduction, and let $T=BA^{-1}$, which is a transvection in all cases considered. Then we prove that one has either $G=\langle T\rangle\ast\langle B\rangle$ or $G=\langle\pm T\rangle\ast_{\{\pm I\}}\langle B\rangle$ if $-I\in\langle B\rangle$.
To see that this implies that $G$ is thin, note that $\langle T\rangle$, $\langle \pm T\rangle$ and $\langle B\rangle$, as abelian groups, have the Haagerup property \cite[Ex.~2.4(2)]{Jolissaint}, whence also $G$ is Haagerup \cite[Prop.~2.5(3)]{Jolissaint}.
But if $G$ were a finite-index subgroup of $\mathrm{Sp}_6(\mathbb{Z})$, then it would be a lattice in $\mathrm{Sp}_6(\mathbb{R})$ and hence satisfy property (T) \cite[Thm.~1.5.3, Thm.~1.7.1]{bhv}. This is a contradiction, since non-compact groups cannot both be Haagerup and have property (T) \cite[Sect.~1.2.1]{CCJJV01}.

To prove that $G=\langle T\rangle\ast\langle B\rangle$ or $G=\langle \pm T\rangle\ast_{\{\pm I\}}\langle B\rangle$, we adapt the strategy of Brav and Thomas~\cite{BT} -- playing ping-pong on a set of convex cones in $\mathbb{R}^n$ -- from dimension $n=4$ to $n=6$.
Namely, we apply the following version of the ping pong lemma to the case where $G_1=\langle T\rangle$, $G_2=\langle B\rangle$, $H=\{I\}$, respectively to $G_1=\langle\pm T\rangle$, $G_2=\langle B\rangle$, $H=\{\pm I\}$.

\begin{thm}[see \cite{LS77}, Prop.~III.12.4]\label{th:pingpong}
Let $G$ be a group generated by two subgroups $G_1$ and $G_2$ whose intersection $H$ has index $>2$ in $G_{1}$ or $G_{2}$. Suppose that $G$ acts on a set $W$, and suppose that there are disjoint non-empty subsets $X$ and $Y$ such that $(G_1\setminus H) Y \subseteq X$ and $(G_2\setminus H) X \subseteq Y$, with $H Y \subseteq Y$ and $H X \subseteq X$. Then $G = G_1 \star_{H} G_2$.
\end{thm}

The two halves $X$, $Y$ of our ping pong table will both decompose into two parts $X=X^+\cup X^-$, $Y=Y^+\cup Y^-$, which in turn decompose into a union of open convex cones.

Let $P_{\pi}$ be the permutation matrix that sends the basis vector $e_i$ to $e_{\pi(i)}$ for $i\in \{1,\dots,6\}$, where $\pi=(1\, 6)(2\, 5)(3\, 4)$.
Then $E=BP_\pi$ satisfies

\begin{align}\label{eq:Eprop}
E^2 & =I, & EBE^{-1} & =B^{-1}, & ETE^{-1} & =T^{-1}.
\end{align}
Let also
\begin{equation*}
\eta=\min\{k>0\mid (B^{k}-I)^{6}=0\}\in\mathbb{N}^{+}
\end{equation*}
and set
\begin{equation}\label{eq:table-halves}
\begin{aligned}
X^+& =C\cup -C, & X^- & =EX^+, \\
Y^+ & =\bigcup\{\pm B^i D \mid 1\leq i\leq\eta, B^i\neq\pm I\}, & Y^- & =EY^+,
\end{aligned}
\end{equation}
for some open convex cones $C$ and $D$. This structure for the ping pong table halves is the same that was used by Brav and Thomas, where it was motivated by the geometric picture of the two-dimensional case, see \cite[Sect.~2]{BT}. Unlike in \cite{BT}, however, the cone $D$ is not necessarily defined starting from $C$: this allows us to cover $6$ more cases.

With this setup we see that the conditions $H Y \subseteq Y$, $H X \subseteq X$ in the ping pong lemma are automatic. In order to satisfy the other hypotheses of Theorem~\ref{th:pingpong}, it is sufficient that the following conditions hold:
\begin{enumerate}
\item\label{it:nonempty} $X\neq\emptyset$,
\item\label{it:disjoint} $X\cap Y=\emptyset$,
\item\label{it:Bcheck} $B^{i}X\subseteq Y$ for all $i\in\mathbb{Z}, B^i\neq\pm I$,
\item\label{it:Tcheck} $T^i Y\subseteq X$ for all $i\in\mathbb{Z},i\neq 0$.
\end{enumerate}
Condition \eqref{it:Tcheck} involves infinitely many checks, but it can be replaced with
\begin{enumerate}
\renewcommand{\theenumi}{4\alph{enumi}}
\renewcommand{\labelenumi}{(\theenumi)}
\item\label{it:Tcheck_plus} $T^{-1}(Y\cup X^{+})\subseteq X^{+}$,
\item\label{it:Tcheck_minus} $T(Y\cup X^{-})\subseteq X^{-}$.
\end{enumerate}
If $B$ has finite order, \eqref{it:Bcheck} involves only finitely many checks. If $B$ has infinite order, we proceed like we did with \eqref{it:Tcheck} and replace \eqref{it:Bcheck} with

\begin{enumerate}
\renewcommand{\theenumi}{3\alph{enumi}}
\renewcommand{\labelenumi}{(\theenumi)}
\item\label{it:Bcheck_plus} $B(X\cup Y^{+})\subseteq Y^{+}$,
\item\label{it:Bcheck_minus} $B^{-1}(X\cup Y^{-})\subseteq Y^{-}$.
\end{enumerate}

Figure~\ref{fig:ppt} shows a visualization of the containment conditions for our ping pong.

It remains to choose the cones $C$ and $D$. We will give a cone $C$ for each individual case in terms of finitely many generating rays. A suitable cone $D$ can then usually be directly obtained from $C$: when $B$ has finite order, we simply set $D=C$, and when $B$ has infinite order, it often suffices to let $D$ be the cone spanned by $C$, $EC$, $(B^\eta-I)C$, $(B^\eta-I)EC$. In the few cases where this does not work, we will also supply a working choice for $D$. A general description of the method yielding $C$ and $D$ is given in Section~\ref{se:construction}, and the actual choice of cones is given in Sections~\ref{se:thinlist}-\ref{se:morethin}-\ref{se:thinsp4}.

\begin{center}
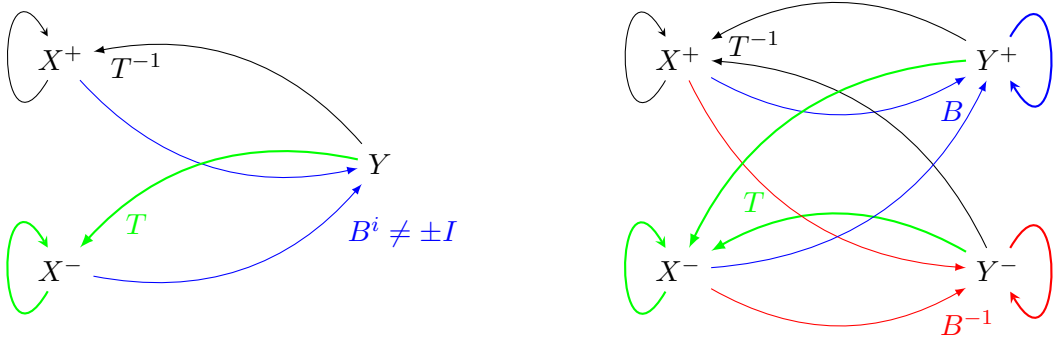
\begin{figure}%
\centering
\begin{minipage}{7cm}
\centering
\begin{tikzpicture}[scale=0.70]
\node (a) at (0,0) {$X^{+}$};
\node (b) at (6,-2) {$Y$};
\draw[blue,-latex,bend right]  (a) edge (b);
\draw[-latex,bend right]  (b) edge (a);
\node (c) at (0,-4) {$X^{-}$};
\draw[blue,-latex,bend right]  (c) edge (b);
\draw[green, thick,-latex,bend right]  (b) edge (c);
\node[xshift=1cm,yshift=-0.1cm] at (a) {$T^{-1}$};
\node[xshift=1cm,yshift=0.6cm] at (c) {\color{green}$T$};
\node[xshift=0.3cm,yshift=-0.9cm] at (b) {\color{blue}$B^{i}\neq\pm I$};
\node[xshift=0.4cm,yshift=-0.7cm] at (b) {\color{blue}};
\draw[->,>= stealth]  (a) edge [out=240,in=120,distance=20mm]   (a);
\draw[green,thick, ->,>= stealth]  (c) edge [out=240,in=120,distance=20mm]   (c);
\end{tikzpicture}\\

\end{minipage}
\hfill
\begin{minipage}{7cm}
\centering
\begin{tikzpicture}[scale=0.70]
\node (a) at (0,0) {$X^{+}$};
\node (b) at (6,0) {$Y^{+}$};
\draw[blue,-latex,bend right]  (a) edge (b);
\draw[-latex,bend right]  (b) edge (a);
\node (c) at (0,-4) {$X^{-}$};
\node (d) at (6,-4) {$Y^{-}$};
\draw[red,-latex,bend right]  (c) edge (d);
\draw[green, thick,-latex,bend right]  (d) edge (c);
\draw[red,-latex,bend right ] (a) edge (d);
\draw[-latex, bend right] (d) edge (a);
\draw[green, thick,-latex,bend right ] (b) edge (c);
\draw[blue,-latex, bend right] (c) edge (b);
\node[xshift=1cm,yshift=0.2cm] at (a) {$T^{-1}$};
\node[xshift=1cm,yshift=0.9cm] at (c) {\color{green}$T$};
\node[xshift=-0.6cm,yshift=-0.7cm] at (b) {\color{blue}$B$};
\node[xshift=-0.4cm,yshift=-0.7cm] at (d) {\color{red}$B^{-1}$};
\draw[->,>= stealth]  (a) edge [out=240,in=120,distance=20mm]   (a);
\draw[green,thick, ->,>= stealth]  (c) edge [out=240,in=120,distance=20mm]   (c);
\draw[blue,thick, ->,>= stealth]  (b) edge [out=240,in=120,distance=-30mm]   (b);
\draw[red,thick, ->,>= stealth]  (d) edge [out=240,in=120,distance=-30mm]   (d);
\end{tikzpicture}\\
\end{minipage}
\caption{Containment conditions for the ping pong lemma: the left figure is for $B$ of finite order, depicting conditions~\eqref{it:Bcheck}, \eqref{it:Tcheck_plus}, \eqref{it:Tcheck_minus}, and the right figure is for $B$ of infinite order, depicting conditions \eqref{it:Bcheck_plus}, \eqref{it:Bcheck_minus}, \eqref{it:Tcheck_plus}, \eqref{it:Tcheck_minus}.}\label{fig}
\label{fig:ppt}
\end{figure}
\end{center}

\subsection{Construction of the cones}\label{se:construction}

In this section we explain how we found good cones $C,D$ from which the ping pong table can be built. All steps were carried out with the computer.
The basic strategy is that we start with some initial cones $C_0,D_0$ that are small enough that the ping pong table resulting from \eqref{eq:table-halves} will satisfy condition \eqref{it:disjoint} (disjointness). Then we iteratively expand the cones until also conditions \eqref{it:Bcheck} and \eqref{it:Tcheck} (closedness under ping pong) are satisfied.

In the $i$-th expansion step, let the current cone $C_i$ be generated by a finite set $R_i$.
In order to enforce conditions \eqref{it:Tcheck_plus} and \eqref{it:Tcheck_minus}, the new cone $C_{i+1}$ must contain a multiple of the vectors $T^{-1}B^{j}v$ and $T^{-1}B^{j}Ev$ for all $v\in R_{i}$.
We let the generating set $R_{i+1}$ for $C_{i+1}$ consist of $R_{i}$ and all the vectors $\varepsilon_{j}T^{-1}B^{j}v$ and $\varepsilon'_{j}T^{-1}B^{j}Ev$, for all $v\in R_{i}$ and all $j\leq\eta$ and for some appropriate signs $\varepsilon_{j},\varepsilon'_{j}$, determined by the sign of the largest eigenvalues of $T^{-1}B^{j}$ and $T^{-1}B^{j}E$. As a result, the table constructed from $C_{i+1},D_i$ now satisfies condition \eqref{it:Tcheck}. Similarly, we pass from $D_i$ to $D_{i+1}$, such that the table obtained from $C_{i+1},D_{i+1}$ satisfies \eqref{it:Bcheck}. The expansion of $D$ may interfere with condition \eqref{it:Tcheck}, and so we continue with expanding $C$ again. We want that the sequences $(C_i)_{i},(D_i)_{i}$ become stationary after finitely many steps and that the ping pong table halves $X,Y$ resulting from the limit cones $C,D$ are still disjoint.

For this, we need to make a good choice for the initial cones $C_0,D_0$.
Consider first the rank-$1$ operator $T-I$. Its image $\mathrm{Im}(T-I)$ intersects non-trivially the closure of any orbit $\{T^{-i}v\mid i\in\mathbb{N}\}$, where $v\notin\mathrm{ker}(T-I)$. Thus, $\mathrm{Im}(T-I)$ must be contained in the closure of $X^{+}$ by~\eqref{it:Tcheck_plus}.
The cones $C_0,D_0$ should both contain a generator $t_0$ of the subspace $\mathrm{Im}(T-I)$ (its sign can be chosen freely at this point).

If the sequences $(C_i)_{i},(D_i)_{i}$ starting with $C_0=D_0=\mathbb{R}_{+}t_0$ (where $\mathbb{R}_{+}$ is the set of non-negative reals) do not appear to become stationary, we add more rays to $C_0,D_0$. The minimal way to do this is to identify boundary points of the old $\bigcup_{i}C_i$ and $\bigcup_{i}D_i$, and add them to $C_0$ and $D_0$, respectively.

One particularly important boundary point is the limit point
\[u_0\vcentcolon=\lim_{i\to\infty}\, (T^{-1}B)^{\lambda i}\, t_0=\lim_{i\to\infty}\, A^{\lambda i}\, t_0,\]
where $\lambda\in\mathbb{N}$ is the smallest number such that $A^\lambda$ is unipotent. In our setup -- with the particular generators $A,B$ of certain hypergeometric groups -- the point $u_0$ must lie in the intersection of the closures of both ping pong table halves.
So, $C$ must not contain a neighborhood of $u_0$, but it must contain the whole set $\{A^{\lambda i} t_0\}_{i\in\mathbb{N}}$, which approaches $u_0$ from a direction determined by the powers $(A^\lambda-I)^k$.
To deal with this, we add to $C_0$ convex combinations $\sum_k\mu_k(A^\lambda-I)^k t_0$, where the weights $\mu_k$ will be heavier for higher $k$. This is similar to how Brav and Thomas \cite{BT} make use of the logarithm.
If $B$ has infinite order, we proceed similarly with the unipotent matrix $B^\eta$ and the cone $D_0$.

\medskip

The cones $C$ and $D$ produced by the method above were often defined by hundreds of rays, which seemed unnecessarily complicated. To reduce the number of rays, we attempted to simplify or remove individual generating rays of $C,D$ in a brute force search. Another approach that was very successful in practice was to enlarge the initial cone $C_0$ by adding a random vector in the vicinity of $t_0$.
Bigger starting cones $C_0,D_0$ often led to very fast convergence of the expansion procedure and to much simpler solutions.
However, if the starting cones are too big, the ping pong table halves will start to overlap at some point during the expansion procedure, violating condition \eqref{it:disjoint} and causing the ping pong to fail.
We note that in some cases this happens even for the minimal choice $C_0=D_0=\mathbb{R}_{+}t_0$, in which case our ping pong approach cannot work; see Section~\ref{se:open}.

\subsection{SageMath code}\label{se:code}

Here we present the SageMath code that performs the verification. Its inputs are the polynomials $f$ and $g$, a matrix whose columns represent the rays spanning the cone $C$, and a potential second matrix for the cone $D$. If the output is ``True'', the resulting objects satisfy the conditions of the ping pong lemma. The run-time is very short: for each case, the code takes only a few seconds.

There are two main functions, for the cases of $B$ of finite and infinite order respectively. They are given below, together with the relevant auxiliary functions. The same code can also be accessed on GitHub~\cite{BDN21-s}.

\hspace{0.3cm}
\small
\begin{verbatim}
from itertools import count
def companion_matrix(polynomial):
    n=polynomial.degree()
    left_block=block_matrix([[matrix(1,n-1)],[matrix.identity(n-1)]])
    return block_matrix([[left_block,-matrix(polynomial.list()[:-1]).T]])

##### Partial checks #####
# In the following, CC and DD are unions of open cones,
# and they are given as a list of such cones.
def are_disjoint(CC, DD): # checks whether CC,DD are disjoint
    return all(not C.intersection(D).is_solid() for C in CC for D in DD)
def contained_in(CC,DD): # checks whether CC is contained in DD
    P=[any([all([D.contains(r) for r in C.rays()]) for D in DD]) for C in CC]
    return all(P)
def transform_set(CC,L): # applies the linear transformation L to CC
    return [Cone((L*C.rays().column_matrix()).T) for C in CC]

##### Main functions: verifying hypotheses of the ping-pong lemma #####
def verify_finite_order(f,g,J):
    B=companion_matrix(g)
    T=B*companion_matrix(f).inverse()
    E=B*Permutation(range(f.degree(),0,-1)).to_matrix()
    eta=next(k for k in count(1) if (B^k-1)^f.degree()==0)
    if not B^eta==1:
        return False # B is not of finite order
    # build ping pong table
    Xplus=[Cone(J.T),Cone(-J.T)]
    Xminus=[Cone((E*J).T),Cone(-(E*J).T)]
    Y=[C for k in [1..eta-1] for C in transform_set(Xplus+Xminus, B^k)
       if not B^k==-1]
    # verify ping pong
    if not Xplus[0].is_solid():
        return False # verifies (1)
    if not are_disjoint(Xplus+Xminus, Y):
        return False # verifies (2)
    if not all(contained_in(transform_set(Xplus+Xminus, B^ii), Y) 
               for ii in [1..eta-1] if not B^ii==-1):
        return False # verifies (3)
    if not contained_in(transform_set(Y+Xplus, T.inverse()), Xplus):
        return False # verifies (4a)
    if not contained_in(transform_set(Y+Xminus, T), Xminus):
        return False # verifies (4b)
    return True

def verify_infinite_order(f,g,J,K=None):
    B=companion_matrix(g)
    T=B*companion_matrix(f).inverse()
    E=B*Permutation(range(f.degree(),0,-1)).to_matrix()
    eta=next(k for k in count(1) if (B^k-1)^f.degree()==0)
    if K==None:
        K=block_matrix([[J,E*J,(B^eta-1)*J,(B^eta-1)*E*J]])
    # build ping pong table
    Xplus=[Cone(J.T),Cone(-J.T)]
    Xminus=[Cone((E*J).T),Cone(-(E*J).T)]
    Ytemp=[Cone(K.T),Cone(-K.T)]
    Yplus=[C for k in [1..eta] for C in transform_set(Ytemp, B^k)]
    Yminus=transform_set(Yplus, E)
    # verify ping pong for infinite order
    if not Xplus[0].is_solid():
        return False # verifies (1)
    if not are_disjoint(Xplus+Xminus, Yplus+Yminus):
        return False # verifies (2)
    if not contained_in(transform_set(Xplus+Xminus+Yplus, B), Yplus):
        return False # verifies (3a)
    if not contained_in(transform_set(Xplus+Xminus+Yminus, B.inverse()), Yminus):
        return False # verifies (3b)
    if not contained_in(transform_set(Yplus+Yminus+Xplus, T.inverse()), Xplus):
        return False # verifies (4a)
    if not contained_in(transform_set(Yplus+Yminus+Xminus, T), Xminus):
        return False # verifies (4b)
    return True
\end{verbatim}
\normalsize

\hspace{0.3cm}

In SageMath, a cone is defined by a matrix whose rows are the vectors spanning the rays of the cone. In our paper, both for convenience and for consistency with the existing literature, the rays of the cones are instead spanned by the columns of our matrices. Furthermore, the cones in SageMath are closed, unlike ours. Consequently, the disjointness of our (open) cones is checked by verifying that the intersection of two (closed) cones is less than $6$-dimensional.
In our code and in the SageMath class ConvexRationalPolyhedralCone that implements the cones, all computations are done with rational numbers. Hence, no rounding errors can occur and the calculations are exact.

To complete the process, one has only to write down the inputs after the definitions above, and evaluate the appropriate main function. We use \verb|verify_finite_order(f,g,M)| when $B$ has finite order, \verb|verify_infinite_order(f,g,M)| when $B$ has infinite order and $D$ is implicitly defined from $C$ as described in Section~\ref{se:pingpong}, and \verb|verify_infinite_order(f,g,M,M1)| when we provide $D$ ourselves. Here are two examples of the additional code necessary to verify thinness for cases A-2 and A-37 of table \ref{tab:thin}.

\small
\begin{verbatim}
### Case A-2 ###
f=cyclotomic_polynomial(1)^6
g=cyclotomic_polynomial(2)^4*cyclotomic_polynomial(3)
M=matrix([
[  0,   0,   0,   0,   1,   1,   1,   4,   5],
[ -1,  -1,  -1,  -1,  -5,  -5,  -1, -19, -19],
[ -1,   0,   0,   4,  10,  10,  10,  42,  32],
[ -4,  -1,  -1,  -6, -10, -10, -13, -58, -32],
[ -1,   0,   1,   4,   5,   7,   7,  46,  19],
[ -1,  -1,  -2,  -1,  -1,  -3,  -4, -15,  -5]])
M1=matrix([
[ -1,  -1,  -1,   0,   0,   0,   0,   0,   1,   1,   1,   4,   5],
[ -5,  -4,  -1,  -2,  -1,  -1,  -1,   0,  -5,  -5,  -1, -19, -19],
[ -8,  -7,  -1,   1,  -3,   0,   4,   0,  10,  10,  10,  42,  32],
[-14,  -7,  -1,  -1,  -4,  -1,  -6,  -2, -10, -10, -13, -58, -32],
[-11,  -4,  -1,   0,  -3,   1,   4,   1,   5,   7,   7,  46,  19],
[ -5,  -1,  -1,  -1,  -1,  -2,  -1,  -1,  -1,  -3,  -4, -15,  -5]])
print(verify_infinite_order(f,g,M,M1))

### Case A-37 ###
f=cyclotomic_polynomial(1)^6
g=cyclotomic_polynomial(7)
M=matrix([
[   0,   0,   0,    1,    4,   25],
[ -40,  -2,  -1,   -5,  -20,  -45],
[  83,   0,   4,   10,   55,  106],
[ -86,  -1,  -6,  -10,  -85, -158],
[  83,   0,   4,    5,   65,  133],
[ -40,  -2,  -1,   -1,  -19,  -61]])
print(verify_finite_order(f,g,M))
\end{verbatim}
\normalsize

\section{Ping pong tables for maximally unipotent \texorpdfstring{$\Sp(6)$}{Sp(6)} hypergeometric groups}\label{se:thinlist}

In this section, we prove Theorem~\ref{thm:total} by establishing the thinness of the $17$ hypergeometric groups with a maximally unipotent monodromy marked as ``thin'' in Table~\ref{tab:thin}. For all these cases, we have $\alpha=\left(0,0,0,0,0,0\right)$. For each case we indicate whether the order of $B$ is finite or infinite, and we write the matrix $M$ whose columns represent the rays defining the cone $C$. When the cone $D$ is not automatically defined as described in Section~\ref{se:pingpong}, we also provide the additional matrix $M_{1}$ defining $D$. The list of matrices is also accessible on GitHub \cite{BDN21-s}, arranged so as to be readily usable for our verification code.

\subsection{Thinness of case A-1}
$B$ has infinite order, and the matrices are

\resizebox{\textwidth}{!}{
\parbox{\textwidth}{
\begin{align*}
M  = &\left(\begin{array}{rrrrrrrrrrrrrrrr}
0 & 0 & 0 & 0 & 1 & 1 & 1 & 1 & 1 & 1 & 1 & 2 & 2 & 3 & 3 & 7 \\
-3 & -3 & -1 & -1 & -16 & -11 & -10 & -6 & -5 & -5 & -2 & -14 & -9 & -18 & -15 & -39 \\
0 & 0 & 0 & 4 & 5 & 11 & 22 & 14 & 10 & 11 & 1 & 27 & 16 & 42 & 42 & 90 \\
-10 & -8 & -6 & -6 & -40 & -36 & -31 & -22 & -10 & -32 & -8 & -58 & -23 & -52 & -50 & -106 \\
0 & 0 & 0 & 4 & 0 & 3 & 20 & 12 & 5 & 5 & 0 & 19 & 15 & 35 & 35 & 63 \\
-4 & -2 & -2 & -1 & -9 & -10 & -9 & -8 & -1 & -13 & -2 & -25 & -5 & -10 & -15 & -15
\end{array}\right),\\
M_{1}  =& \left(\begin{array}{rrrrrrrrrrrrrrrr}
-61 & -45 & -11 & -11 & -10 & -9 & -7 & -6 & -6 & -4 & -3 & -1 & -1 & -1 & -1 & -1 \\
-355 & -256 & -65 & -60 & -59 & -51 & -65 & -34 & -33 & -21 & -17 & -14 & -8 & -7 & -6 & -5 \\
-850 & -593 & -158 & -132 & -144 & -126 & -159 & -73 & -83 & -44 & -39 & -3 & -15 & -10 & -5 & -10 \\
-1062 & -718 & -210 & -146 & -186 & -154 & -184 & -79 & -105 & -46 & -46 & -42 & -26 & -30 & -27 & -10 \\
-705 & -450 & -135 & -81 & -126 & -105 & -105 & -41 & -70 & -24 & -28 & -1 & -14 & -5 & -6 & -5 \\
-231 & -120 & -61 & -18 & -51 & -35 & -27 & -9 & -22 & -5 & -7 & -14 & -7 & -11 & -9 & -1
\end{array}\right.
 & & & & & & & & & & & & & \\
& \left. \begin{array}{rrrrrrrrrrrrrrrr}
0 & 0 & 0 & 0 & 0 & 0 & 0 & 1 & 1 & 1 & 1 & 2 & 2 & 3 & 3 & 7 \\
-4 & -3 & -2 & -2 & -1 & -1 & -1 & -6 & -5 & -2 & 0 & -14 & -9 & -18 & -15 & -39 \\
7 & 0 & -2 & 0 & 0 & 2 & 4 & 14 & 10 & 1 & 1 & 27 & 16 & 42 & 42 & 90 \\
-11 & -10 & -4 & -4 & -6 & -5 & -6 & -22 & -10 & -8 & -7 & -58 & -23 & -52 & -50 & -106 \\
7 & 0 & -2 & 0 & 0 & 1 & 4 & 12 & 5 & 0 & 0 & 19 & 15 & 35 & 35 & 63 \\
-3 & -4 & -1 & -1 & -2 & -2 & -1 & -8 & -1 & -2 & -3 & -25 & -5 & -10 & -15 & -15
\end{array}\right).
\end{align*}
}}

\subsection{Thinness of case A-2}
$B$ has infinite order, and the matrices are
\begin{align*}
M & =\left(\begin{array}{rrrrrrrrr}
0 & 0 & 0 & 0 & 1 & 1 & 1 & 4 & 5 \\
-1 & -1 & -1 & -1 & -5 & -5 & -1 & -19 & -19 \\
-1 & 0 & 0 & 4 & 10 & 10 & 10 & 42 & 32 \\
-4 & -1 & -1 & -6 & -10 & -10 & -13 & -58 & -32 \\
-1 & 0 & 1 & 4 & 5 & 7 & 7 & 46 & 19 \\
-1 & -1 & -2 & -1 & -1 & -3 & -4 & -15 & -5
\end{array}\right),\\
M_{1} & =\left(\begin{array}{rrrrrrrrrrrrr}
-1 & -1 & -1 & 0 & 0 & 0 & 0 & 0 & 1 & 1 & 1 & 4 & 5 \\
-5 & -4 & -1 & -2 & -1 & -1 & -1 & 0 & -5 & -5 & -1 & -19 & -19 \\
-8 & -7 & -1 & 1 & -3 & 0 & 4 & 0 & 10 & 10 & 10 & 42 & 32 \\
-14 & -7 & -1 & -1 & -4 & -1 & -6 & -2 & -10 & -10 & -13 & -58 & -32 \\
-11 & -4 & -1 & 0 & -3 & 1 & 4 & 1 & 5 & 7 & 7 & 46 & 19 \\
-5 & -1 & -1 & -1 & -1 & -2 & -1 & -1 & -1 & -3 & -4 & -15 & -5
\end{array}\right).
\end{align*}

\subsection{Thinness of case A-3}
$B$ has infinite order, and the matrix is
\begin{align*}
M & =\left(\begin{array}{rrrrrrrrr}
0 & 0 & 0 & 0 & 1 & 1 & 1 & 1 & 1 \\
-1 & -1 & -1 & -1 & -6 & -5 & -4 & -4 & -2 \\
-2 & 0 & 4 & 4 & 13 & 10 & 7 & 8 & 1 \\
-2 & -1 & -7 & -6 & -14 & -10 & -7 & -8 & 0 \\
-2 & 0 & 6 & 4 & 8 & 5 & 4 & 5 & 0 \\
-1 & -1 & -2 & -1 & -2 & -1 & -1 & -2 & -4
\end{array}\right).
\end{align*}

\subsection{Thinness of case A-4}
$B$ has infinite order, and the matrix is
\begin{align*}
M & =\left(\begin{array}{rrrrrr}
0 & 0 & 0 & 1 & 1 & 1 \\
-1 & -1 & -1 & -5 & -3 & -1 \\
-1 & 2 & 4 & 10 & 4 & 3 \\
0 & -2 & -6 & -10 & -4 & -5 \\
-1 & 2 & 4 & 5 & 3 & 4 \\
-1 & -1 & -1 & -1 & -1 & -2
\end{array}\right).
\end{align*}

\subsection{Thinness of case A-5}
$B$ has infinite order, and the matrix is
\[
\resizebox{\textwidth}{!}{
$M=\left(\begin{array}{rrrrrrrrrrrrr}
0 & 0 & 0 & 0 & 0 & 0 & 1 & 1 & 1 & 1 & 5 & 11 & 23 \\
-1 & -1 & -1 & -1 & -1 & -1 & -6 & -5 & -3 & -1 & -25 & -66 & -127 \\
0 & 0 & 0 & 1 & 2 & 4 & 14 & 10 & 8 & 0 & 66 & 154 & 290 \\
-3 & -3 & -2 & -3 & -3 & -6 & -18 & -10 & -14 & -6 & -81 & -188 & -338 \\
0 & 1 & 0 & 2 & 1 & 4 & 12 & 5 & 9 & 1 & 56 & 123 & 199 \\
-1 & -2 & -1 & -1 & -1 & -1 & -5 & -1 & -6 & -3 & -21 & -34 & -47
\end{array}\right)$.
}
\]
\subsection{Thinness of case A-6}
$B$ has infinite order, and the matrix is
$$
M=\left(\begin{array}{rrrrrrr}
0 & 0 & 0 & 0 & 1 & 1 & 7 \\
-6 & -1 & -1 & -1 & -5 & 0 & -25 \\
17 & 0 & 0 & 4 & 10 & 0 & 50 \\
-22 & -2 & -2 & -6 & -10 & 0 & -70 \\
17 & 0 & 0 & 4 & 5 & 0 & 55 \\
-6 & -2 & -1 & -1 & -1 & -1 & -17
\end{array}\right).
$$

\subsection{Thinness of case A-7}
$B$ has infinite order, and the matrix is
$$ M=\left(\begin{array}{rrrrrr}
0 & 0 & 0 & 1 & 1 & 7 \\
-1 & -1 & -1 & -5 & 0 & -25 \\
0 & 2 & 4 & 10 & 0 & 50 \\
-1 & -2 & -6 & -10 & -1 & -70 \\
0 & 2 & 4 & 5 & 2 & 55 \\
-1 & -1 & -1 & -1 & -2 & -17
\end{array}\right).
$$

\subsection{Thinness of case A-8}
$B$ has infinite order, and the matrix is
$$ M=\left(\begin{array}{rrrrrrr}
0 & 0 & 0 & 1 & 1 & 1 & 12 \\
-2 & -1 & -1 & -5 & -1 & -1 & -33 \\
0 & 3 & 4 & 10 & 1 & 4 & 44 \\
-3 & -4 & -6 & -10 & -1 & -5 & -54 \\
0 & 3 & 4 & 5 & 1 & 3 & 48 \\
-2 & -1 & -1 & -1 & -1 & -2 & -17
\end{array}\right).$$

\subsection{Thinness of case A-9}
$B$ has infinite order, and the matrix is
$$ M=\left(\begin{array}{rrrrrrr}
0 & 0 & 0 & 0 & 1 & 1 & 1 \\
-1 & -1 & -1 & -1 & -5 & -4 & -1 \\
0 & 0 & 2 & 4 & 10 & 8 & 2 \\
-1 & -1 & -2 & -6 & -10 & -10 & -2 \\
0 & 1 & 2 & 4 & 5 & 7 & 1 \\
-1 & -1 & -1 & -1 & -1 & -2 & -1
\end{array}\right).
$$

\subsection{Thinness of case A-10}
$B$ has infinite order, and the matrix is
$$ M=\left(\begin{array}{rrrrrrrr}
0 & 0 & 0 & 0 & 0 & 1 & 1 & 277 \\
-8 & -1 & -1 & -1 & -1 & -5 & -2 & -685 \\
25 & 0 & 0 & 0 & 4 & 10 & 5 & 970 \\
-34 & -2 & -1 & 0 & -6 & -10 & -6 & -1570 \\
25 & 0 & 0 & 0 & 4 & 5 & 4 & 1585 \\
-8 & -1 & -2 & -1 & -1 & -1 & -2 & -577
\end{array}\right).
$$

\subsection{Thinness of case A-11}
$B$ has infinite order, and the matrix is
$$M=\left(\begin{array}{rrrrrr}
0 & 0 & 0 & 1 & 1 & 21 \\
-1 & -1 & -1 & -5 & -2 & -78 \\
1 & 3 & 4 & 10 & 4 & 131 \\
-2 & -4 & -6 & -10 & -5 & -135 \\
1 & 3 & 4 & 5 & 4 & 84 \\
-1 & -1 & -1 & -1 & -2 & -23
\end{array}\right).
$$

\subsection{Thinness of case A-12}
$B$ has infinite order, and the matrix is
$$M=\left(\begin{array}{rrrrrrr}
0 & 0 & 0 & 0 & 1 & 1 & 3 \\
-1 & -1 & -1 & -1 & -5 & 0 & -10 \\
0 & 0 & 2 & 4 & 10 & 0 & 20 \\
-1 & 0 & -2 & -6 & -10 & 0 & -30 \\
0 & 0 & 2 & 4 & 5 & 0 & 25 \\
-2 & -1 & -1 & -1 & -1 & -1 & -8
\end{array}\right).
$$

\subsection{Thinness of case A-13}
$B$ has infinite order, and the matrix is
$$M=\left(\begin{array}{rrrrrrr}
0 & 0 & 0 & 0 & 1 & 1 & 4 \\
-1 & -1 & -1 & -1 & -5 & -1 & -15 \\
0 & 0 & 2 & 4 & 10 & 1 & 30 \\
0 & 0 & -2 & -6 & -10 & -1 & -40 \\
0 & 1 & 2 & 4 & 5 & 1 & 30 \\
-1 & -1 & -1 & -1 & -1 & -1 & -9
\end{array}\right).
$$

\subsection{Thinness of case A-14}
$B$ has infinite order, and the matrix is
$$M=\left(\begin{array}{rrrrrrr}
0 & 0 & 0 & 0 & 1 & 1 & 3 \\
-1 & -1 & -1 & -1 & -5 & 0 & -10 \\
0 & 0 & 2 & 4 & 10 & 0 & 20 \\
0 & 0 & -2 & -6 & -10 & 0 & -30 \\
0 & 1 & 2 & 4 & 5 & 0 & 25 \\
-1 & -1 & -1 & -1 & -1 & -1 & -8
\end{array}\right).
$$

\subsection{Thinness of case A-31}
$B$ has finite order, and the matrix is
$$M=\left(\begin{array}{rrrrrr}
0 & 0 & 0 & 1 & 1 & 4 \\
-35 & -4 & -1 & -5 & -4 & -20 \\
78 & 6 & 4 & 10 & 8 & 55 \\
-86 & -7 & -6 & -10 & -8 & -85 \\
78 & 6 & 4 & 5 & 5 & 65 \\
-35 & -4 & -1 & -1 & -2 & -19
\end{array}\right).$$

\subsection{Thinness of case A-37}
$B$ has finite order, and the matrix is
$$
M=\left(\begin{array}{rrrrrr}
0 & 0 & 0 & 1 & 1 & 4 \\
-31 & -1 & -1 & -5 & -1 & -20 \\
74 & 0 & 4 & 10 & 3 & 55 \\
-86 & -1 & -6 & -10 & -5 & -85 \\
74 & 0 & 4 & 5 & 4 & 65 \\
-31 & -1 & -1 & -1 & -2 & -19
\end{array}\right).
$$

\subsection{Thinness of case A-38} 
$B$ has finite order, and the matrix is
$$ M=\left(\begin{array}{rrrrrr}
0 & 0 & 0 & 1 & 1 & 4 \\
-29 & -2 & -1 & -5 & -2 & -20 \\
72 & 2 & 4 & 10 & 5 & 55 \\
-86 & -3 & -6 & -10 & -7 & -85 \\
72 & 2 & 4 & 5 & 5 & 65 \\
-29 & -2 & -1 & -1 & -2 & -19
\end{array}\right).$$

\section{Ping pong tables for more \texorpdfstring{$\Sp(6)$}{Sp(6)} hypergeometric groups}\label{se:morethin}

In this section, we prove Theorem~\ref{thm:total-nonmum} by establishing the thinness of the $46$ hypergeometric groups from Table~\ref{tab:thin-more}. The solutions are represented by matrices, in the same way as in Section~\ref{se:thinlist}. The list of matrices is also accessible on GitHub \cite{BDN21-s}.

\subsection{Thinness of case C-2}
$B$ has infinite order, and the matrices are

\resizebox{\textwidth}{!}{
\parbox{\textwidth}{
\begin{align*}
& M=\left(\begin{array}{rrrrrrrrrrrrrrrrr}
0 & 0 & 0 & 1 & 1 & 1 & 1 & 1 & 2 & 3 & 3 & 5 & 5 & 5 & 6 & 13 & 21 \\
-7 & -3 & -1 & -8 & -5 & -3 & -2 & -1 & -17 & -15 & -11 & -28 & -16 & -14 & -37 & -44 & -68 \\
-4 & 0 & 1 & -4 & -2 & 2 & 1 & 1 & -8 & -3 & -1 & -11 & 11 & 12 & -16 & 20 & 42 \\
-10 & -2 & 0 & -11 & -8 & -4 & -1 & -1 & -24 & -20 & -17 & -44 & -16 & -8 & -54 & -57 & -74 \\
-4 & 0 & 1 & -3 & -3 & 1 & 2 & 1 & -8 & -6 & -5 & -16 & 6 & 11 & -14 & 5 & 21 \\
-7 & -3 & -1 & -7 & -7 & -5 & -1 & -1 & -16 & -19 & -16 & -34 & -22 & -14 & -37 & -65 & -94
\end{array}\right),\\
& M_1= \left(\begin{array}{rrrrrrrrrrrrrrrrrrrrrrrrr}
-39 & -30 & -26 & -21 & -13 & -7 & -6 & -5 & -5 & -2 & -1 & -1 & -1 & 0 & 0 & 1 & 1 & 1 & 3 & 5 & 5 & 6 & 13 & 21 \\
-138 & -106 & -91 & -178 & -117 & -24 & -61 & -42 & -34 & -24 & -5 & -5 & -3 & -1 & -1 & -3 & -2 & -1 & -11 & -16 & -14 & -37 & -44 & -68 \\
-207 & -159 & -135 & -126 & -86 & -36 & -56 & -29 & -24 & -22 & -6 & -5 & -4 & -2 & 1 & 2 & 1 & 1 & -1 & 11 & 12 & -16 & 20 & 42 \\
-209 & -161 & -136 & -242 & -161 & -36 & -102 & -56 & -48 & -40 & -9 & -9 & -4 & -2 & 0 & -4 & -1 & -1 & -17 & -16 & -8 & -54 & -57 & -74 \\
-166 & -128 & -109 & -105 & -71 & -29 & -58 & -24 & -23 & -22 & -6 & -6 & -3 & -2 & 1 & 1 & 2 & 1 & -5 & 6 & 11 & -14 & 5 & 21 \\
-73 & -56 & -47 & -152 & -96 & -12 & -61 & -36 & -34 & -25 & -5 & -6 & -1 & -1 & -1 & -5 & -1 & -1 & -16 & -22 & -14 & -37 & -65 & -94
\end{array}\right).
\end{align*}
}}

\subsection{Thinness of case C-3}
$B$ has infinite order, and the matrices are

\resizebox{\textwidth}{!}{
\parbox{\textwidth}{
\begin{align*}
& M=\left(\begin{array}{rrrrrrrrrrr}
0 & 0 & 0 & 0 & 1 & 1 & 1 & 1 & 1 & 1 & 3 \\
-3 & -2 & -1 & -1 & -7 & -6 & -5 & -3 & -2 & -1 & -9 \\
0 & 1 & 0 & 1 & 0 & 1 & 0 & 0 & 1 & 1 & 8 \\
-2 & -2 & -3 & 0 & -5 & -3 & -3 & -1 & -1 & -1 & -6 \\
0 & 1 & 1 & 1 & 1 & 0 & 1 & 0 & 2 & 1 & 6 \\
-3 & -2 & -3 & -1 & -6 & -5 & -6 & -3 & -1 & -1 & -8
\end{array}\right),\\
& M_1= \left(\begin{array}{rrrrrrrrrrrrrrr}
-87 & -75 & -54 & -14 & -3 & -1 & -1 & -1 & 0 & 0 & 0 & 0 & 1 & 1 & 3 \\
-210 & -183 & -129 & -33 & -17 & -4 & -4 & -2 & -3 & -1 & -1 & -1 & -2 & -1 & -9 \\
-135 & -120 & -81 & -21 & -3 & -2 & -1 & -1 & 1 & -1 & 0 & 1 & 1 & 1 & 8 \\
-89 & -78 & -54 & -14 & -12 & -3 & -3 & -1 & -3 & 0 & -3 & 0 & -1 & -1 & -6 \\
-206 & -180 & -127 & -33 & -1 & -2 & -2 & -2 & 0 & -1 & 1 & 1 & 2 & 1 & 6 \\
-137 & -122 & -83 & -21 & -18 & -4 & -5 & -1 & -1 & -1 & -3 & -1 & -1 & -1 & -8
\end{array}\right).
\end{align*}
}}

\subsection{Thinness of case C-4}
$B$ has infinite order, and the matrix is

\resizebox{\textwidth}{!}{
\parbox{\textwidth}{
\begin{align*}
& M =\left(\begin{array}{rrrrrrrrrrrrrrr}
0 & 0 & 0 & 0 & 0 & 0 & 1 & 1 & 1 & 1 & 2 & 3 & 16 & 48 & 49 \\
-62 & -61 & -2 & -2 & -1 & -1 & -3 & -2 & -1 & 0 & -46 & -28 & -93 & -114 & -109 \\
58 & 1 & 0 & 0 & 0 & 1 & 1 & 1 & 0 & 1 & 14 & 6 & 57 & 37 & 47 \\
-63 & -49 & -2 & -2 & -2 & 0 & -3 & -1 & -2 & -1 & -65 & -23 & -65 & -39 & -40 \\
13 & 70 & 0 & 0 & 0 & 1 & -1 & 2 & 1 & 0 & 55 & 0 & 30 & 65 & 46 \\
-49 & -64 & -3 & -1 & -1 & -1 & -3 & -1 & -1 & -1 & -63 & -13 & -48 & -100 & -96
\end{array}\right).
\end{align*}
}}
\subsection{Thinness of case C-5}
$B$ has infinite order, and the matrix is
\begin{align*}
& M=\left(\begin{array}{rrrrrrrr}
0 & 0 & 0 & 0 & 0 & 0 & 1 & 1 \\
-2 & -1 & -1 & -1 & -1 & -1 & -2 & 0 \\
1 & -1 & 0 & 0 & 0 & 1 & 1 & 1 \\
-2 & -1 & -2 & -1 & 0 & 0 & -1 & -1 \\
0 & 0 & 1 & -1 & 0 & 1 & 2 & 0 \\
-1 & -1 & -2 & -1 & -1 & -1 & -1 & -1
\end{array}\right).
\end{align*}

\subsection{Thinness of case C-6}
$B$ has infinite order, and the matrix is
\begin{align*}
& M= \left(\begin{array}{rrrrrrrrr}
0 & 0 & 0 & 0 & 0 & 0 & 0 & 1 & 1 \\
-2 & -2 & -2 & -1 & -1 & -1 & -1 & -2 & 0 \\
1 & 2 & 2 & 0 & 0 & 1 & 1 & 1 & 1 \\
0 & -1 & 0 & -1 & 0 & -1 & 0 & -1 & -1 \\
2 & 2 & 1 & 1 & 0 & 0 & 1 & 2 & 0 \\
-2 & -2 & -2 & -1 & -1 & -1 & -1 & -1 & -1
\end{array}\right).
\end{align*}

\subsection{Thinness of case C-7}
$B$ has infinite order, and the matrix is

\resizebox{\textwidth}{!}{
\parbox{\textwidth}{
\begin{align*}
& M=\left(\begin{array}{rrrrrrrrrrrrrrr}
0 & 0 & 0 & 0 & 0 & 1 & 1 & 3 & 69 & 78 & 82 & 91 & 91 & 140 & 140 \\
-195 & -2 & -1 & -1 & -1 & -2 & -1 & -140 & -173 & -42 & -404 & -195 & -133 & -419 & -322 \\
280 & 0 & 0 & 1 & 1 & 1 & 1 & 223 & 0 & 0 & 280 & 231 & 29 & 224 & 0 \\
-302 & 0 & 0 & -1 & 0 & -1 & -1 & -305 & -20 & -182 & -33 & -182 & 49 & 0 & 0 \\
244 & 1 & 0 & 1 & 1 & 2 & 1 & 280 & 313 & 182 & 102 & 91 & 42 & 117 & 371 \\
-140 & -1 & -1 & -2 & -1 & -1 & -1 & -174 & -302 & -91 & -140 & -91 & -133 & -175 & -302
\end{array}\right).
\end{align*}
}}
\subsection{Thinness of case C-8}
$B$ has infinite order, and the matrix is

\begin{align*}
& M=\left(\begin{array}{rrrrrrrrrr}
0 & 0 & 0 & 0 & 0 & 1 & 1 & 1 & 1 & 2 \\
-2 & -1 & -1 & -1 & -1 & -2 & -2 & -2 & -1 & -2 \\
1 & 0 & 1 & 1 & 1 & 1 & 1 & 3 & 1 & 1 \\
0 & 0 & -1 & 0 & 0 & -1 & 0 & -2 & -1 & -1 \\
1 & 1 & 1 & 0 & 1 & 2 & 1 & 2 & 2 & 1 \\
-1 & -1 & -1 & -1 & -1 & -1 & -2 & -2 & -2 & -1
\end{array}\right).
\end{align*}

\subsection{Thinness of case C-11}
$B$ has infinite order, and the matrices are

\resizebox{\textwidth}{!}{
\parbox{\textwidth}{
\begin{align*}
& M=\left(\begin{array}{rrrrrrrrrrrrrrr}
0 & 0 & 0 & 0 & 0 & 1 & 1 & 1 & 1 & 1 & 1 & 3 & 3 & 6 \\
-9 & -4 & -4 & -2 & -1 & -6 & -3 & -3 & -3 & -1 & -1 & -18 & -11 & -21 \\
-4 & 0 & 0 & 0 & 2 & -1 & -1 & 4 & 6 & 1 & 3 & 0 & 18 & 31 \\
-22 & -10 & -7 & -5 & -2 & -12 & -8 & -4 & -7 & -1 & -4 & -61 & -18 & -30 \\
-4 & 0 & -1 & 0 & 2 & -2 & -1 & 3 & 6 & 1 & 0 & -11 & 14 & 24 \\
-9 & -5 & -3 & -2 & -1 & -4 & -3 & -1 & -3 & -1 & -4 & -29 & -6 & -10
\end{array}\right),\\
& M_1=\left(\begin{array}{rrrrrrrrrrrrrrrrrrr}
-1994 & -415 & -6 & -3 & -1 & -1 & -1 & -1 & 0 & 0 & 0 & 0 & 0 & 1 & 1 & 1 & 1 & 3 & 6 \\
-8911 & -1859 & -40 & -21 & -8 & -6 & -6 & -4 & -5 & -4 & -3 & -1 & -1 & -3 & -3 & -1 & -1 & -11 & -21 \\
-17117 & -3580 & -42 & -19 & -5 & -10 & -8 & -7 & -17 & 0 & -1 & -3 & 2 & 4 & 6 & 1 & 3 & 18 & 31 \\
-18485 & -3876 & -114 & -60 & -21 & -15 & -18 & -7 & -24 & -7 & -7 & -4 & -2 & -4 & -7 & -1 & -4 & -18 & -30 \\
-11550 & -2429 & -35 & -15 & -5 & -10 & -7 & -4 & -19 & -1 & 0 & -3 & 2 & 3 & 6 & 1 & 0 & 14 & 24 \\
-3335 & -705 & -51 & -26 & -8 & -6 & -8 & -1 & -7 & -3 & -4 & -1 & -1 & -1 & -3 & -1 & -4 & -6 & -10
\end{array}\right).
\end{align*}
}}

\subsection{Thinness of case C-12}
$B$ has infinite order, and the matrix is

\resizebox{\textwidth}{!}{
\parbox{\textwidth}{
\begin{align*}
& M =\left(\begin{array}{rrrrrrrrrrr}
0 & 0 & 0 & 0 & 0 & 0 & 1 & 1 & 77 & 165 & 198 \\
-321 & -1 & -1 & -1 & -1 & -1 & -3 & 0 & -308 & -583 & -683 \\
52 & -1 & 0 & 0 & 1 & 2 & 4 & 1 & -462 & 847 & 887 \\
0 & -1 & -3 & -1 & -5 & -2 & -4 & -1 & -473 & -1782 & -835 \\
52 & 0 & 0 & -1 & 0 & 2 & 3 & 0 & 448 & 682 & -298 \\
-666 & -1 & -1 & -1 & -2 & -1 & -1 & -1 & -165 & -212 & -152
\end{array}\right).
\end{align*}
}}

\subsection{Thinness of case C-13}
$B$ has infinite order, and the matrix is

\resizebox{\textwidth}{!}{
\parbox{\textwidth}{
\begin{align*}
& M=\left(\begin{array}{rrrrrrrrrr}
0 & 0 & 0 & 0 & 0 & 1 & 1 & 1 & 3 & 6 \\
-2 & -1 & -1 & -1 & -1 & -3 & -3 & -1 & -11 & -21 \\
2 & 0 & 0 & 1 & 2 & 4 & 6 & 1 & 18 & 31 \\
-3 & -2 & -1 & -2 & -2 & -4 & -7 & -1 & -18 & -30 \\
1 & 1 & 0 & 1 & 2 & 3 & 6 & 1 & 14 & 24 \\
-2 & -2 & -1 & -1 & -1 & -1 & -3 & -1 & -6 & -10
\end{array}\right).
\end{align*}
}}
\subsection{Thinness of case C-14}
$B$ has infinite order, and the matrix is

\resizebox{\textwidth}{!}{
\parbox{\textwidth}{
\begin{align*}
& M=\left(\begin{array}{rrrrrrrrrrrr}
0 & 0 & 0 & 0 & 1 & 1 & 1 & 1 & 1 & 3 & 3 & 6 \\
-1 & -1 & -1 & -1 & -6 & -5 & -3 & -3 & -1 & -12 & -11 & -21 \\
0 & 0 & 1 & 2 & 0 & 8 & 4 & 6 & 1 & 18 & 18 & 31 \\
-2 & 0 & -1 & -2 & -11 & -8 & -4 & -7 & -1 & -30 & -18 & -30 \\
0 & 0 & 1 & 2 & 9 & 0 & 3 & 6 & 1 & 14 & 14 & 24 \\
-1 & -1 & -5 & -1 & -4 & -3 & -1 & -3 & -1 & -6 & -6 & -10
\end{array}\right).
\end{align*}
}}

\subsection{Thinness of case C-15}
$B$ has infinite order, and the matrix is

\resizebox{\textwidth}{!}{
\parbox{\textwidth}{
\begin{align*}
M= & \left(\begin{array}{rrrrrrrrrrrrrrrrrrrrrrrrrrrrrrrrrr}
0 & 0 & 0 & 0 & 0 & 1 & 1 & 2 & 2 & 4 & 7 & 10 & 12 & 12 & 14 & 17 & 36 & 41 & 63 & 76 & 132 & 151 \\
-5 & -3 & -2 & -1 & -1 & -3 & -2 & -10 & -8 & -14 & -24 & -51 & -65 & -58 & -49 & -80 & -127 & -152 & -216 & -263 & -465 & -528 \\
11 & 6 & 3 & 2 & 2 & 4 & 3 & 15 & 14 & 18 & 35 & 74 & 116 & 86 & 68 & 126 & 172 & 222 & 314 & 380 & 708 & 794 \\
-13 & -10 & -5 & -3 & -2 & -4 & -3 & -16 & -15 & -17 & -38 & -76 & -134 & -88 & -69 & -132 & -162 & -212 & -332 & -386 & -742 & -828 \\
9 & 7 & 3 & 2 & 2 & 3 & 2 & 9 & 10 & 12 & 32 & 48 & 88 & 53 & 54 & 81 & 120 & 153 & 279 & 320 & 592 & 671 \\
-5 & -4 & -2 & -1 & -1 & -1 & -1 & -4 & -6 & -7 & -16 & -29 & -41 & -30 & -30 & -36 & -63 & -76 & -132 & -151 & -249 & -284
\end{array}\right).
\end{align*}
}}

\subsection{Thinness of case C-16}
$B$ has infinite order, and the matrix is

\resizebox{\textwidth}{!}{
\parbox{\textwidth}{
\begin{align*}
& M=\left(\begin{array}{rrrrrrrrrrr}
0 & 0 & 0 & 0 & 1 & 1 & 1 & 1 & 2 & 3 & 6 \\
-4 & -2 & -2 & -1 & -10 & -3 & -3 & -1 & -12 & -11 & -21 \\
0 & 1 & 1 & 2 & 20 & 4 & 6 & 1 & 20 & 18 & 31 \\
-8 & -2 & -1 & -2 & -26 & -4 & -7 & -1 & -20 & -18 & -30 \\
10 & 1 & 0 & 2 & 20 & 3 & 6 & 1 & 14 & 14 & 24 \\
-11 & -2 & -1 & -1 & -10 & -1 & -3 & -1 & -7 & -6 & -10
\end{array}\right).
\end{align*}
}}

\subsection{Thinness of case C-17}
$B$ has infinite order, and the matrix is

\resizebox{\textwidth}{!}{
\parbox{\textwidth}{
\begin{align*}
&M=\left(\begin{array}{rrrrrrrrrrrrr}
0 & 0 & 0 & 0 & 1 & 1 & 1 & 1 & 1 & 1 & 3 & 3 & 6 \\
-7 & -1 & -1 & -1 & -7 & -5 & -3 & -3 & -1 & -1 & -11 & -11 & -21 \\
8 & 1 & 1 & 2 & 12 & 6 & 4 & 6 & 1 & 2 & 16 & 18 & 31 \\
-5 & -2 & -1 & -2 & -12 & -7 & -4 & -7 & -1 & -2 & -18 & -18 & -30 \\
4 & 1 & 1 & 2 & 7 & 7 & 3 & 6 & 1 & 1 & 14 & 14 & 24 \\
-3 & -1 & -1 & -1 & -3 & -3 & -1 & -3 & -1 & -1 & -5 & -6 & -10
\end{array}\right).
\end{align*}
}}

\subsection{Thinness of case C-18}
$B$ has infinite order, and the matrix is

\resizebox{\textwidth}{!}{
\parbox{\textwidth}{
\begin{align*}
&M=\left(\begin{array}{rrrrrrrrrrrrrrrrr}
0 & 0 & 0 & 0 & 0 & 0 & 0 & 0 & 1 & 1 & 1 & 2 & 3 & 6 & 7 & 12 & 14 \\
-6 & -6 & -5 & -2 & -2 & -1 & -1 & -1 & -3 & -3 & -1 & -9 & -11 & -21 & -26 & -41 & -44 \\
5 & 9 & 7 & 2 & 4 & 0 & 0 & 2 & 4 & 6 & 1 & 12 & 18 & 31 & 40 & 58 & 57 \\
-5 & -14 & -5 & -1 & -5 & -1 & -1 & -2 & -4 & -7 & -1 & -15 & -18 & -30 & -44 & -56 & -54 \\
7 & 10 & 7 & 1 & 3 & 0 & 3 & 2 & 3 & 6 & 1 & 16 & 14 & 24 & 34 & 40 & 42 \\
-5 & -5 & -5 & -1 & -2 & -1 & -3 & -1 & -1 & -3 & -1 & -7 & -6 & -10 & -12 & -14 & -16
\end{array}\right).
\end{align*}
}}

\subsection{Thinness of case C-19}
$B$ has finite order, and the matrix is
$$ M=\left(\begin{array}{rrrrrrrrrrrrr}
0 & 0 & 0 & 0 & 0 & 1 & 1 & 1 & 1 & 2 & 2 & 3 & 10 \\
-2 & -1 & -1 & -1 & -1 & -3 & -3 & -1 & -1 & -9 & -9 & -11 & -37 \\
1 & 0 & 0 & 0 & 2 & 4 & 6 & 1 & 5 & 12 & 12 & 18 & 59 \\
-1 & -2 & -1 & -1 & -2 & -4 & -3 & -3 & -5 & -16 & -16 & -18 & -62 \\
0 & 0 & 0 & 1 & 2 & 3 & 2 & 2 & 1 & 8 & 12 & 18 & 52 \\
-1 & -1 & -1 & -1 & -1 & -1 & -3 & -1 & -1 & -4 & -5 & -10 & -22
\end{array}\right).$$

\subsection{Thinness of case C-20}
$B$ has finite order, and the matrix is
\begin{align*}
& M=\left(\begin{array}{rrrrrrrrrrrrr}
0 & 0 & 0 & 0 & 1 & 1 & 1 & 1 & 1 & 1 & 2 & 3 & 6 \\
-2 & -1 & -1 & -1 & -5 & -5 & -5 & -3 & -3 & -2 & -7 & -10 & -21 \\
4 & 1 & 2 & 2 & 7 & 8 & 9 & 3 & 4 & 3 & 9 & 14 & 32 \\
-6 & -1 & -3 & -2 & -6 & -8 & -10 & -3 & -4 & -3 & -8 & -15 & -34 \\
6 & 1 & 2 & 2 & 4 & 5 & 7 & 3 & 3 & 2 & 6 & 13 & 27 \\
-3 & -1 & -1 & -1 & -2 & -2 & -3 & -2 & -1 & -1 & -3 & -6 & -11
\end{array}\right).
\end{align*}

\subsection{Thinness of case C-21}
$B$ has infinite order, and the matrices are

\resizebox{\textwidth}{!}{
\parbox{\textwidth}{
\begin{align*}
& M=\left(\begin{array}{rrrrrrrrrrrrrrrr}
0 & 0 & 0 & 0 & 1 & 1 & 1 & 1 & 1 & 1 & 8 & 33 & 94 & 209 \\
-4 & -3 & -1 & -1 & -7 & -7 & -5 & -4 & -2 & -1 & -39 & -157 & -437 & -951 \\
1 & 0 & 0 & 3 & 4 & 14 & 10 & 7 & 0 & 1 & 83 & 324 & 877 & 1862 \\
-7 & -7 & -4 & -4 & -25 & -23 & -17 & -7 & -4 & -1 & -102 & -379 & -992 & -2049 \\
0 & 0 & 0 & 3 & 1 & 12 & 7 & 4 & 0 & 1 & 71 & 261 & 655 & 1307 \\
-2 & -3 & -1 & -1 & -10 & -10 & -8 & -1 & -2 & -1 & -33 & -94 & -209 & -390
\end{array}\right),\\
&M_1= \left(\begin{array}{rrrrrrrrrrrrrrrrrrrrrrr}
-599 & -390 & -209 & -94 & -34 & -33 & -8 & -2 & -1 & -1 & -1 & 0 & 0 & 0 & 0 & 1 & 1 & 1 & 1 & 8 & 33 & 94 & 209 \\
-2692 & -1741 & -1435 & -679 & -151 & -259 & -73 & -29 & -6 & -6 & -4 & -4 & -2 & -1 & -1 & -7 & -5 & -4 & -1 & -39 & -157 & -437 & -951 \\
-5201 & -3339 & -992 & -379 & -289 & -102 & -17 & 3 & -10 & -7 & -7 & 1 & 0 & -3 & 3 & 14 & 10 & 7 & 1 & 83 & 324 & 877 & 1862 \\
-5647 & -3598 & -4975 & -2308 & -310 & -841 & -214 & -74 & -15 & -21 & -7 & -7 & -7 & -4 & -4 & -23 & -17 & -7 & -1 & -102 & -379 & -992 & -2049 \\
-3548 & -2241 & -437 & -157 & -193 & -39 & -5 & 7 & -10 & -4 & -4 & 0 & 1 & -3 & 3 & 12 & 7 & 4 & 1 & 71 & 261 & 655 & 1307 \\
-1033 & -643 & -1996 & -907 & -55 & -322 & -79 & -25 & -6 & -9 & -1 & -2 & -4 & -1 & -1 & -10 & -8 & -1 & -1 & -33 & -94 & -209 & -390
\end{array}\right).
\end{align*}
}}
\subsection{Thinness of case C-22}
$B$ has infinite order, and the matrix is

\resizebox{\textwidth}{!}{
\parbox{\textwidth}{
\begin{align*}
& M=\left(\begin{array}{rrrrrrrrrrrrrrrrrr}
0 & 0 & 0 & 0 & 0 & 0 & 0 & 1 & 1 & 1 & 1 & 2 & 7 & 28 & 68 & 123 & 189 \\
-7 & -2 & -2 & -1 & -1 & -1 & -1 & -4 & -4 & -1 & -1 & -8 & -35 & -133 & -312 & -547 & -822 \\
7 & 1 & 4 & 0 & 0 & 0 & 3 & 7 & 8 & 0 & 0 & 20 & 70 & 273 & 615 & 1041 & 1532 \\
-5 & -6 & -5 & -3 & -2 & -1 & -4 & -7 & -9 & -6 & -2 & -26 & -91 & -322 & -679 & -1107 & -1605 \\
7 & 0 & 3 & 0 & 1 & 0 & 3 & 4 & 6 & 0 & 0 & 19 & 72 & 217 & 426 & 674 & 972 \\
-7 & -2 & -2 & -1 & -1 & -1 & -1 & -1 & -2 & -4 & -2 & -9 & -28 & -68 & -123 & -189 & -271
\end{array}\right).
\end{align*}
}}
\subsection{Thinness of case C-23}
$B$ has infinite order, and the matrix is

\resizebox{\textwidth}{!}{
\parbox{\textwidth}{
\begin{align*}
& M= \left(\begin{array}{rrrrrrrrrrrrrrrrrrr}
0 & 0 & 0 & 0 & 1 & 1 & 2 & 2 & 14 & 15 & 46 & 53 & 61 & 143 & 176 & 358 & 402 \\
-1 & -1 & -1 & -1 & -4 & -1 & -12 & -11 & -68 & -73 & -227 & -251 & -290 & -662 & -819 & -1707 & -1834 \\
0 & 0 & 1 & 3 & 7 & 1 & 17 & 20 & 143 & 153 & 448 & 515 & 598 & 1322 & 1646 & 3316 & 3603 \\
-2 & -2 & -3 & -4 & -7 & -1 & -36 & -32 & -176 & -193 & -492 & -599 & -701 & -1487 & -1866 & -3594 & -3982 \\
0 & 1 & 1 & 3 & 4 & 1 & 14 & 17 & 122 & 129 & 311 & 407 & 478 & 974 & 1235 & 2242 & 2556 \\
-1 & -2 & -1 & -1 & -1 & -1 & -15 & -14 & -53 & -61 & -89 & -143 & -176 & -308 & -402 & -630 & -775
\end{array}\right).
\end{align*}
}}

\subsection{Thinness of case C-24}
$B$ has infinite order, and the matrix is

\resizebox{\textwidth}{!}{
\parbox{\textwidth}{
\begin{align*}
& M=\left(\begin{array}{rrrrrrrrrrrrrrr}
0 & 0 & 0 & 0 & 0 & 0 & 1 & 1 & 1 & 1 & 3 & 8 & 8 \\
-5 & -1 & -1 & -1 & -1 & -1 & -5 & -4 & -1 & -1 & -18 & -44 & -41 \\
12 & 0 & 0 & 1 & 1 & 3 & 10 & 7 & 1 & 1 & 38 & 89 & 80 \\
-14 & -3 & -2 & -3 & -1 & -4 & -13 & -7 & -1 & -1 & -45 & -100 & -86 \\
9 & 0 & 1 & 0 & 1 & 3 & 10 & 4 & 1 & 1 & 30 & 62 & 51 \\
-3 & -1 & -1 & -1 & -1 & -1 & -4 & -1 & -2 & -1 & -9 & -17 & -14
\end{array}\right).
\end{align*}
}}
\subsection{Thinness of case C-25}
$B$ has infinite order, and the matrix is

\resizebox{\textwidth}{!}{
\parbox{\textwidth}{
\begin{align*}
& M=\left(\begin{array}{rrrrrrrrrrrrr}
0 & 0 & 0 & 0 & 0 & 1 & 1 & 1 & 1 & 5 & 15 & 33 & 58 \\
-1 & -1 & -1 & -1 & -1 & -4 & -4 & -1 & -1 & -24 & -70 & -150 & -265 \\
0 & 0 & 0 & 0 & 3 & 7 & 10 & 1 & 1 & 51 & 141 & 293 & 508 \\
-2 & -2 & -2 & -1 & -4 & -7 & -13 & -1 & -1 & -60 & -159 & -321 & -545 \\
0 & 0 & 1 & 0 & 3 & 4 & 10 & 0 & 1 & 42 & 105 & 204 & 337 \\
-2 & -1 & -2 & -1 & -1 & -1 & -5 & -1 & -1 & -15 & -33 & -60 & -94
\end{array}\right).
\end{align*}
}}

\subsection{Thinness of case C-26}
$B$ has infinite order, and the matrix is

\resizebox{\textwidth}{!}{
\parbox{\textwidth}{
\begin{align*}
& M= \left(\begin{array}{rrrrrrrrrrrrrrrr}
0 & 0 & 0 & 0 & 0 & 0 & 1 & 1 & 1 & 1 & 5 & 13 & 24 & 34 \\
-5 & -2 & -2 & -1 & -1 & -1 & -6 & -4 & -1 & -1 & -24 & -60 & -107 & -158 \\
10 & 1 & 1 & 0 & 0 & 3 & 10 & 7 & 1 & 1 & 49 & 119 & 204 & 297 \\
-12 & -4 & -1 & -3 & 0 & -4 & -14 & -7 & -1 & -1 & -60 & -133 & -217 & -311 \\
10 & 0 & 0 & 2 & 0 & 3 & 12 & 4 & 1 & 1 & 41 & 83 & 131 & 187 \\
-5 & -2 & -1 & -2 & -1 & -1 & -5 & -1 & -2 & -1 & -13 & -24 & -37 & -51
\end{array}\right).
\end{align*}
}}

\subsection{Thinness of case C-27}
$B$ has infinite order, and the matrix is

\resizebox{\textwidth}{!}{
\parbox{\textwidth}{
\begin{align*}
& M= \left(\begin{array}{rrrrrrrrrrrr}
0 & 0 & 0 & 0 & 1 & 1 & 1 & 4 & 8 & 9 & 14 \\
-5 & -1 & -1 & -1 & -6 & -4 & -1 & -19 & -35 & -41 & -61 \\
11 & 0 & 1 & 3 & 10 & 7 & 1 & 38 & 63 & 80 & 113 \\
-13 & -1 & -1 & -4 & -13 & -7 & -1 & -46 & -62 & -88 & -116 \\
11 & 0 & 1 & 3 & 11 & 4 & 1 & 31 & 35 & 53 & 66 \\
-5 & -1 & -1 & -1 & -4 & -1 & -1 & -9 & -10 & -14 & -17
\end{array}\right).
\end{align*}
}}

\subsection{Thinness of case C-28}
$B$ has infinite order, and the matrix is

\resizebox{\textwidth}{!}{
\parbox{\textwidth}{
\begin{align*}
& M= \left(\begin{array}{rrrrrrrrrrrrrr}
0 & 0 & 0 & 0 & 0 & 1 & 1 & 1 & 2 & 2 & 9 & 21 & 33 & 40 \\
-6 & -2 & -1 & -1 & -1 & -7 & -4 & -1 & -11 & -9 & -43 & -96 & -144 & -167 \\
11 & 3 & 0 & 1 & 3 & 8 & 7 & 1 & 20 & 18 & 88 & 188 & 267 & 296 \\
-11 & -2 & -1 & -2 & -4 & -10 & -7 & -1 & -27 & -23 & -106 & -206 & -274 & -293 \\
11 & 1 & 0 & 3 & 3 & 10 & 4 & 1 & 24 & 17 & 72 & 125 & 157 & 166 \\
-6 & -1 & -1 & -2 & -1 & -6 & -1 & -1 & -9 & -6 & -21 & -33 & -40 & -43
\end{array}\right).
\end{align*}
}}

\subsection{Thinness of case C-33}
$B$ has finite order, and the matrix is
$$M=\left(\begin{array}{rrrrrrrrrrrr}
0 & 0 & 0 & 0 & 1 & 1 & 1 & 1 & 1 & 2 & 9 & 24 \\
-4 & -2 & -1 & -1 & -4 & -4 & -3 & -3 & -1 & -9 & -43 & -111 \\
9 & 1 & 1 & 3 & 7 & 10 & 2 & 8 & 2 & 21 & 90 & 221 \\
-13 & -1 & -3 & -4 & -7 & -12 & -4 & -11 & -1 & -26 & -105 & -246 \\
10 & 1 & 1 & 3 & 4 & 8 & 3 & 8 & 1 & 21 & 73 & 159 \\
-4 & -2 & -1 & -1 & -1 & -3 & -2 & -3 & -2 & -9 & -24 & -47
\end{array}\right). $$

\subsection{Thinness of case C-34}
$B$ has finite order, and the matrix is

\resizebox{\textwidth}{!}{
\parbox{\textwidth}{
\begin{align*}
& M=\left(\begin{array}{rrrrrrrrrrrrrr}
0 & 0 & 0 & 0 & 1 & 1 & 1 & 1 & 1 & 2 & 4 & 8 & 12 & 15 \\
-5 & -2 & -1 & -1 & -6 & -6 & -6 & -4 & -1 & -17 & -19 & -36 & -52 & -63 \\
8 & 2 & 2 & 3 & 10 & 12 & 13 & 7 & 1 & 36 & 38 & 69 & 96 & 113 \\
-9 & -3 & -3 & -4 & -13 & -13 & -16 & -7 & -1 & -44 & -43 & -74 & -99 & -114 \\
4 & 2 & 2 & 3 & 9 & 9 & 12 & 4 & 1 & 29 & 28 & 45 & 58 & 66 \\
-2 & -2 & -1 & -1 & -4 & -3 & -4 & -1 & -1 & -9 & -8 & -12 & -15 & -17
\end{array}\right).
\end{align*}
}}

\subsection{Thinness of case C-35}
$B$ has infinite order, and the matrix is

\resizebox{\textwidth}{!}{
\parbox{\textwidth}{
\begin{align*}
 M= &\left(\begin{array}{rrrrrrrrrrrrrrrrrrrrrrrrrrrrrrrr}
0 & 0 & 0 & 0 & 1 & 1 & 1 & 2 & 2 & 2 & 7 & 12 & 53 & 84 & 84 & 91 & 105 & 129 & 140 & 169 & 187 & 189 & 284 & 288 & 304 \\
-2 & -1 & -1 & -1 & -7 & -1 & 0 & -2 & 0 & 0 & -5 & -6 & -140 & 0 & 56 & 15 & -21 & 155 & -120 & 119 & -410 & -84 & 20 & -99 & -239 \\
0 & -1 & -1 & 0 & 0 & -1 & 1 & -1 & 0 & 1 & 0 & 0 & 208 & 140 & -36 & 210 & 105 & 149 & 0 & 70 & -328 & 168 & 45 & 204 & 319 \\
-2 & -2 & -1 & -2 & -11 & -2 & -1 & -6 & -4 & -6 & -16 & -35 & -173 & -204 & -168 & -175 & -70 & -213 & -105 & -134 & -635 & -273 & -249 & -408 & -424 \\
0 & -1 & -1 & 0 & -1 & -1 & 0 & -1 & -3 & -6 & 0 & -7 & 56 & -84 & -21 & 20 & -99 & -120 & 0 & -239 & 30 & 119 & -140 & 15 & 155 \\
-1 & -1 & -1 & -2 & -15 & -1 & -1 & -2 & -2 & -4 & -7 & -11 & -173 & -105 & -84 & -330 & -189 & -269 & -84 & -154 & -393 & -288 & -129 & -169 & -284
\end{array}\right).
\end{align*}
}}

\subsection{Thinness of case C-36}
$B$ has infinite order, and the matrix is

\resizebox{\textwidth}{!}{
\parbox{\textwidth}{
\begin{align*}
& M= \left(\begin{array}{rrrrrrrrrrrrrrrrrrrrrrrrrr}
0 & 0 & 0 & 1 & 1 & 1 & 1 & 1 & 1 & 1 & 1 & 1 & 1 & 1 & 2 & 2 & 2 & 2 & 3 & 3 & 3 & 3 \\
-1 & -1 & -1 & -2 & -2 & -1 & -1 & -1 & 0 & 0 & 0 & 0 & 0 & 0 & -1 & 0 & 1 & 1 & -2 & -2 & -1 & 0 \\
-1 & -1 & 0 & -1 & -1 & 0 & 0 & 2 & 0 & 0 & 1 & 1 & 1 & 1 & 2 & 3 & 1 & 2 & 1 & 3 & 2 & 1 \\
-3 & -1 & -1 & -5 & -4 & -2 & -1 & -2 & -2 & -2 & -2 & -2 & -1 & -1 & -3 & -3 & -2 & -3 & -7 & -4 & -4 & -3 \\
-1 & -1 & 0 & 0 & 0 & 0 & 0 & 0 & -1 & 0 & -1 & -1 & -1 & 0 & 1 & 0 & -2 & -1 & -2 & 1 & 0 & -1 \\
-1 & -1 & -1 & -3 & -3 & -1 & -1 & -2 & -1 & -1 & -2 & -1 & -2 & -1 & -3 & -4 & -2 & -3 & -3 & -3 & -2 & -2
\end{array}\right).
\end{align*}
}}

\subsection{Thinness of case C-37}
$B$ has infinite order, and the matrix is

\begin{align*}
&M=\left(\begin{array}{rrrrrrrrrrr}
0 & 0 & 0 & 0 & 0 & 0 & 1 & 1 & 1 & 1 \\
-1 & -1 & -1 & -1 & -1 & -1 & -2 & -1 & -1 & -1 \\
-1 & 0 & 0 & 0 & 0 & 1 & 1 & 0 & 1 & 1 \\
-1 & -2 & -1 & -1 & 0 & -1 & -2 & -1 & -2 & -1 \\
-1 & 0 & 0 & 1 & 0 & 0 & 1 & 0 & 1 & 1 \\
-1 & -1 & -2 & -1 & -1 & -1 & -1 & -1 & -2 & -1
\end{array}\right).
\end{align*}

\subsection{Thinness of case C-38}
$B$ has infinite order, and the matrix is

\resizebox{\textwidth}{!}{
\parbox{\textwidth}{
\begin{align*}
 M= &\left(\begin{array}{rrrrrrrrrrrrrrrrrrrrrrrrrrrrrrrrr}
0 & 0 & 0 & 0 & 0 & 1 & 1 & 1 & 1 & 1 & 1 & 1 & 1 & 2 & 2 & 2 & 3 & 3 \\
-8 & -3 & -3 & -3 & -1 & -8 & -4 & -2 & -2 & -2 & -1 & -1 & -1 & -5 & 0 & 0 & -12 & -8 \\
0 & 0 & 1 & 2 & 0 & 0 & 2 & -1 & 1 & 2 & 1 & 1 & 2 & 2 & 0 & 0 & 3 & 0 \\
-9 & -1 & 0 & -1 & -1 & -4 & -2 & -2 & -4 & -1 & -1 & -1 & -2 & -2 & -3 & -2 & -4 & -8 \\
5 & 0 & 1 & 0 & 0 & 4 & 4 & 1 & 0 & 1 & 0 & 1 & 0 & 4 & 0 & -1 & 9 & 6 \\
-6 & -1 & -3 & -3 & -2 & -5 & -4 & -2 & -2 & -2 & -3 & -1 & -2 & -4 & -2 & -2 & -10 & -6
\end{array}\right.\\
&\left.\begin{array}{rrrrrrrrrrrrrrrrrrrrrrrrrrrrrrrrr}
4 & 4 & 4 & 4 & 4 & 6 & 7 & 10 & 11 & 14 & 15 & 22 & 26 & 30 & 30 \\
-7 & -6 & -4 & -4 & -3 & -16 & -7 & -17 & -12 & -22 & 0 & -37 & -30 & -34 & -30 \\
4 & 3 & 1 & 2 & 1 & 6 & 22 & 8 & 5 & 12 & 0 & 37 & 15 & 30 & 26 \\
-6 & -6 & -4 & -5 & -3 & -10 & -22 & -17 & -14 & -22 & -4 & -22 & -15 & -45 & -30 \\
6 & 6 & 2 & 2 & 0 & 15 & 18 & 16 & 5 & 18 & 0 & 22 & 30 & 45 & 15 \\
-4 & -4 & -2 & -2 & -3 & -14 & -22 & -11 & -6 & -13 & -15 & -26 & -30 & -30 & -15
\end{array}\right).
\end{align*}
}}

\subsection{Thinness of case C-40}
$B$ has infinite order, and the matrix is

\resizebox{\textwidth}{!}{
\parbox{\textwidth}{
\begin{align*}
&M=\left(\begin{array}{rrrrrrrrrrrrrrrrr}
0 & 0 & 0 & 0 & 0 & 0 & 1 & 6 & 6 & 15 & 17 & 19 & 23 & 30 & 36 \\
-3 & -2 & -2 & -1 & -1 & -1 & -4 & 3 & 7 & -6 & -28 & -2 & -40 & -43 & -49 \\
0 & -3 & 0 & -1 & 0 & 0 & 0 & 12 & 16 & -5 & 54 & 8 & 76 & 62 & 68 \\
-8 & -4 & -4 & -3 & -2 & -2 & -10 & -29 & -16 & -14 & -56 & -8 & -76 & -66 & -68 \\
-4 & -3 & 1 & -1 & 0 & 0 & -2 & 4 & 10 & -5 & 22 & -11 & 53 & 34 & 32 \\
-4 & -2 & -2 & -1 & -2 & -1 & -3 & -17 & -23 & -6 & -30 & -6 & -36 & -38 & -19
\end{array}\right).
\end{align*}
}}
\subsection{Thinness of case C-41}
$B$ has infinite order, and the matrix is

\resizebox{\textwidth}{!}{
\parbox{\textwidth}{
\begin{align*}
&M=\left(\begin{array}{rrrrrrrrrrrrrrrrrrr}
0 & 0 & 0 & 0 & 0 & 1 & 2 & 3 & 3 & 5 & 5 & 6 & 7 & 8 & 8 & 10 & 11 & 13 & 15 \\
-4 & -1 & -1 & -1 & -1 & -7 & 3 & -9 & -9 & -8 & -2 & -9 & -4 & -13 & -5 & -15 & -17 & -23 & -19 \\
0 & -1 & 0 & 0 & 0 & 4 & 2 & 2 & 7 & 18 & 10 & 10 & 6 & 15 & -1 & 22 & 31 & 30 & 28 \\
-5 & -2 & -3 & -2 & -2 & -9 & -2 & -13 & -16 & -18 & -21 & -21 & -6 & -30 & -8 & -22 & -34 & -45 & -29 \\
0 & 0 & 0 & -1 & 0 & -3 & 0 & 5 & 8 & 13 & 7 & 0 & -1 & 11 & -1 & 12 & 12 & 23 & 11 \\
-4 & -1 & -1 & -1 & -2 & -6 & -5 & -8 & -13 & -10 & -11 & -6 & -2 & -11 & -5 & -7 & -15 & -18 & -18
\end{array}\right).
\end{align*}
}}
\subsection{Thinness of case C-43}
$B$ has infinite order, and the matrix is

\resizebox{\textwidth}{!}{
\parbox{\textwidth}{
\begin{align*}
& M=\left(\begin{array}{rrrrrrrrrrrrrrrrrrrrr}
0 & 0 & 0 & 0 & 0 & 1 & 1 & 1 & 1 & 2 & 2 & 7 & 15 & 24 & 53 & 84 & 97 \\
-3 & -3 & -2 & -1 & -1 & -4 & -3 & -2 & -1 & -5 & -4 & -6 & -7 & -65 & -135 & -199 & -219 \\
-1 & 0 & 2 & -1 & 1 & 3 & 4 & 4 & 1 & 8 & 5 & 28 & -3 & 114 & 200 & 285 & 310 \\
-8 & -10 & -6 & -2 & -2 & -12 & -5 & -5 & -1 & -8 & -4 & -45 & -16 & -116 & -204 & -304 & -327 \\
-1 & -1 & 1 & -1 & 1 & 0 & 3 & 4 & 1 & 5 & 2 & 19 & -3 & 75 & 149 & 216 & 205 \\
-3 & -4 & -4 & -1 & -1 & -7 & -3 & -2 & -1 & -2 & -1 & -24 & -7 & -53 & -84 & -103 & -87
\end{array}\right).
\end{align*}
}}

\subsection{Thinness of case C-44}
$B$ has infinite order, and the matrix is

\resizebox{\textwidth}{!}{
\parbox{\textwidth}{
\begin{align*}
&M=\left(\begin{array}{rrrrrrrrrrrrrrr}
0 & 0 & 0 & 0 & 0 & 1 & 1 & 1 & 1 & 2 & 2 & 2 \\
-1 & -1 & -1 & -1 & -1 & -3 & -3 & -2 & -1 & -5 & -5 & -4 \\
0 & 0 & 0 & 0 & 1 & 4 & 4 & 4 & 1 & 7 & 8 & 5 \\
-4 & -2 & -1 & 0 & -1 & -6 & -4 & -5 & -1 & -8 & -8 & -4 \\
1 & 0 & 1 & 0 & 0 & 4 & 1 & 4 & 1 & 4 & 5 & 2 \\
-4 & -1 & -1 & -1 & -1 & -2 & -1 & -2 & -1 & -2 & -2 & -1
\end{array}\right).
\end{align*}
}}
\subsection{Thinness of case C-45}
$B$ has infinite order, and the matrix is

\resizebox{\textwidth}{!}{
\parbox{\textwidth}{
\begin{align*}
 M= & \left(\begin{array}{rrrrrrrrrrrrrrrrrrrrrrrrrrrrrrrrrrrrrrrrrrr}
0 & 0 & 0 & 0 & 0 & 1 & 1 & 2 & 2 & 8 & 11 & 22 & 38 & 49 & 55 & 59 & 60 & 66 & 88 & 108 & 114 & 114 \\
-9 & -5 & -4 & -3 & -1 & -2 & -1 & -5 & -4 & -27 & -55 & -77 & -59 & -142 & -77 & -69 & -169 & -176 & -198 & -210 & -304 & -282 \\
14 & 3 & 0 & 4 & 0 & 4 & 1 & 8 & 5 & 40 & 55 & 83 & 113 & 284 & 77 & 85 & 245 & 253 & 264 & 258 & 511 & 401 \\
-18 & -3 & -6 & -5 & 0 & -5 & -1 & -8 & -4 & -59 & -66 & -121 & -151 & -344 & -66 & -96 & -305 & -313 & -275 & -247 & -571 & -439 \\
6 & 2 & 7 & 3 & 0 & 4 & 1 & 5 & 2 & 39 & 66 & 110 & 124 & 257 & 0 & 48 & 234 & 209 & 127 & 101 & 419 & 265 \\
-5 & -3 & -6 & -3 & -1 & -2 & -1 & -2 & -1 & -18 & -60 & -66 & -114 & -153 & -38 & -76 & -114 & -88 & -55 & -59 & -218 & -108
\end{array}\right).
\end{align*}
}}

\subsection{Thinness of case C-46}
$B$ has finite order, and the matrix is

\resizebox{\textwidth}{!}{
\parbox{\textwidth}{
\begin{align*}
M=\left(\begin{array}{rrrrrrrrrrrrrrrrrrrrr}
0 & 0 & 0 & 0 & 0 & 0 & 1 & 1 & 2 & 3 & 3 & 6 & 6 & 7 & 7 & 10 & 10 & 11 & 11 & 14 & 14 \\
-5 & -3 & -2 & -1 & -1 & -1 & -4 & -3 & -5 & -3 & -2 & -15 & -8 & -18 & -7 & -24 & -19 & -23 & -22 & -35 & -28 \\
7 & 1 & 1 & 0 & 1 & 1 & 5 & 4 & 7 & 7 & 8 & 28 & 11 & 32 & 7 & 35 & 28 & 31 & 32 & 52 & 35 \\
-8 & -4 & -2 & -1 & -2 & -2 & -6 & -6 & -8 & -7 & -11 & -28 & -8 & -35 & -7 & -32 & -28 & -31 & -35 & -52 & -32 \\
5 & 1 & 1 & 1 & 1 & 1 & 4 & 4 & 4 & 7 & 8 & 19 & 2 & 28 & 3 & 22 & 15 & 23 & 24 & 35 & 18 \\
-3 & -2 & -1 & -1 & -2 & -1 & -2 & -2 & -2 & -7 & -6 & -10 & -3 & -14 & -3 & -11 & -6 & -11 & -10 & -14 & -7
\end{array}\right).
\end{align*}
}}

\subsection{Thinness of case C-48}
$B$ has infinite order, and the matrix is

\resizebox{\textwidth}{!}{
\parbox{\textwidth}{
\begin{align*}
& M=\left(\begin{array}{rrrrrrrrrrrrrrrrrrrrr}
0 & 0 & 0 & 0 & 0 & 0 & 0 & 1 & 1 & 2 & 3 & 3 & 3 & 5 & 5 & 8 & 21 & 40 & 58 \\
-8 & -3 & -2 & -1 & -1 & -1 & -1 & -3 & -1 & -5 & -11 & -7 & -2 & -17 & -15 & -30 & -76 & -139 & -192 \\
16 & 3 & 0 & 0 & 0 & 0 & 1 & 7 & 1 & 14 & 21 & 9 & 0 & 29 & 23 & 59 & 138 & 244 & 325 \\
-21 & -6 & -2 & -3 & -3 & -3 & -3 & -9 & -1 & -20 & -23 & -7 & -5 & -29 & -21 & -66 & -151 & -262 & -336 \\
15 & 4 & 0 & 0 & 0 & 1 & 0 & 7 & 1 & 11 & 15 & 3 & 0 & 17 & 11 & 44 & 102 & 169 & 202 \\
-8 & -2 & -1 & -5 & -1 & -1 & -1 & -3 & -1 & -8 & -5 & -1 & -2 & -5 & -3 & -21 & -40 & -58 & -63
\end{array}\right).
\end{align*}
}}

\subsection{Thinness of case C-49}
$B$ has infinite order, and the matrix is

\resizebox{\textwidth}{!}{
\parbox{\textwidth}{
\begin{align*}
&M=\left(\begin{array}{rrrrrrrrrrrrrrrr}
0 & 0 & 0 & 0 & 1 & 1 & 3 & 3 & 5 & 5 & 14 & 119 & 442 & 918 & 1220 & 1284 \\
-6 & -2 & -1 & -1 & -3 & -1 & -11 & -7 & -17 & -15 & -14 & -462 & -1649 & -3230 & -3758 & -4218 \\
6 & 0 & 0 & 1 & 7 & 1 & 21 & 9 & 29 & 23 & 42 & 938 & 3074 & 5695 & 5920 & 7042 \\
-16 & -3 & -2 & -2 & -9 & -1 & -23 & -7 & -29 & -21 & -102 & -1148 & -3482 & -6106 & -5644 & -7145 \\
6 & 1 & 0 & 1 & 7 & 1 & 15 & 3 & 17 & 11 & 34 & 850 & 2388 & 3862 & 2993 & 4166 \\
-5 & -2 & -1 & -1 & -3 & -1 & -5 & -1 & -5 & -3 & -119 & -442 & -918 & -1284 & -876 & -1274
\end{array}\right).
\end{align*}
}}

\subsection{Thinness of case C-50}
$B$ has infinite order, and the matrix is

\resizebox{\textwidth}{!}{
\parbox{\textwidth}{
\begin{align*}
M= &\left(\begin{array}{rrrrrrrrrrrrrrrrrrrrrrrrrrrrrrrrrrrrrrr}
0 & 0 & 0 & 0 & 0 & 0 & 1 & 1 & 1 & 2 & 2 & 2 & 3 & 3 & 5 & 5 & 5 & 7 & 9 & 11 & 11 & 14 & 20 & 23 & 27 \\
-7 & -2 & -2 & -1 & -1 & -1 & -6 & -3 & -1 & -15 & -9 & -8 & -11 & -7 & -19 & -17 & -15 & -23 & -34 & -39 & -30 & -45 & -71 & -65 & -88 \\
8 & 3 & 3 & 0 & 0 & 2 & 7 & 7 & 1 & 28 & 15 & 14 & 21 & 9 & 34 & 29 & 23 & 33 & 63 & 69 & 43 & 73 & 126 & 96 & 145 \\
-8 & -5 & -4 & -1 & 0 & -3 & -10 & -9 & -1 & -31 & -19 & -17 & -23 & -7 & -43 & -29 & -21 & -28 & -75 & -76 & -37 & -71 & -137 & -85 & -144 \\
1 & 4 & 3 & 1 & 0 & 2 & 9 & 7 & 1 & 21 & 16 & 11 & 15 & 3 & 30 & 17 & 11 & 13 & 53 & 45 & 17 & 36 & 85 & 40 & 79 \\
-2 & -2 & -2 & -2 & -1 & -1 & -5 & -3 & -1 & -9 & -9 & -4 & -5 & -1 & -11 & -5 & -3 & -6 & -20 & -14 & -8 & -11 & -27 & -13 & -23
\end{array}\right).
\end{align*}
}}

\subsection{Thinness of case C-52}
$B$ has finite order, and the matrix is

\resizebox{\textwidth}{!}{
\parbox{\textwidth}{
\begin{align*}
M= \left(\begin{array}{rrrrrrrrrrrrrrr}
0 & 0 & 0 & 0 & 1 & 1 & 1 & 2 & 3 & 3 & 3 & 5 & 5 & 5 & 5 \\
-4 & -3 & -2 & -1 & -4 & -3 & -1 & -6 & -11 & -11 & -7 & -17 & -17 & -15 & -15 \\
6 & 4 & 3 & 1 & 7 & 7 & 1 & 8 & 20 & 21 & 9 & 29 & 29 & 23 & 23 \\
-7 & -8 & -5 & -1 & -9 & -9 & -1 & -6 & -23 & -23 & -7 & -30 & -29 & -21 & -21 \\
3 & 6 & 3 & 1 & 7 & 7 & 1 & 2 & 15 & 15 & 3 & 17 & 17 & 10 & 11 \\
-2 & -4 & -2 & -1 & -3 & -3 & -1 & -1 & -5 & -5 & -1 & -5 & -5 & -3 & -3
\end{array}\right).
\end{align*}
}}
\subsection{Thinness of case C-53}
$B$ has finite order, and the matrix is

\resizebox{\textwidth}{!}{
\parbox{\textwidth}{
\begin{align*}
M= & \left(\begin{array}{rrrrrrrrrrrrrrrrrrrrrrrrrrrr}
0 & 0 & 0 & 0 & 1 & 1 & 1 & 2 & 2 & 2 & 3 & 3 & 3 & 3  \\
-3 & -2 & -1 & -1 & -5 & -4 & -2 & -11 & -7 & -4 & -19 & -14 & -11 & -7 \\
6 & 5 & 1 & 2 & 9 & 7 & 4 & 18 & 14 & 5 & 35 & 21 & 20 & 9 \\
-7 & -7 & -1 & -3 & -10 & -9 & -5 & -21 & -16 & -4 & -40 & -23 & -23 & -7 \\
4 & 5 & 1 & 2 & 6 & 7 & 4 & 13 & 11 & 2 & 27 & 13 & 15 & 3  \\
-2 & -2 & -1 & -1 & -2 & -3 & -2 & -5 & -4 & -1 & -10 & -5 & -5 & -2 \\
\end{array}\right.
 & & & & & & & & & & & & & \\
& \left. \begin{array}{rrrrrrrrrrrrrrrrrrrrrrrrrrrr}
4 & 4 & 5 & 5 & 5 & 5 & 5 & 7 & 9 & 10 & 12 & 13 & 17 & 22 \\
-14 & -11 & -17 & -17 & -16 & -15 & -14 & -31 & -24 & -37 & -35 & -42 & -61 & -71 \\
25 & 16 & 26 & 29 & 26 & 23 & 20 & 58 & 37 & 61 & 54 & 67 & 105 & 115 \\
-26 & -14 & -29 & -30 & -25 & -21 & -18 & -67 & -36 & -65 & -53 & -69 & -112 & -115 \\
16 & 7 & 17 & 17 & 14 & 10 & 9 & 46 & 19 & 40 & 27 & 39 & 69 & 64 \\
 -5 & -2 & -7 & -5 & -4 & -3 & -3 & -17 & -9 & -13 & -9 & -12 & -22 & -19
\end{array} \right).
\end{align*}
}}

\subsection{Thinness of case C-54}
$B$ has infinite order, and the matrix is

\resizebox{\textwidth}{!}{
\parbox{\textwidth}{
\begin{align*}
&M=\left(\begin{array}{rrrrrrrrrrrrrrrrrr}
0 & 0 & 0 & 0 & 0 & 0 & 0 & 1 & 1 & 1 & 1 & 1 & 1 & 2 & 2 & 2 & 2 \\
-2 & -2 & -1 & -1 & -1 & -1 & -1 & -2 & -1 & -1 & -1 & -1 & 0 & -3 & -2 & -2 & -2 \\
-1 & -1 & -1 & -1 & 0 & 0 & 1 & 0 & -1 & 0 & 0 & 1 & -1 & 1 & -1 & 0 & 0 \\
-2 & -1 & -1 & 0 & -1 & 0 & -1 & 1 & 0 & 0 & 1 & -1 & 0 & 0 & 1 & -1 & 0 \\
-2 & -1 & 0 & 0 & -1 & 0 & 1 & 0 & 2 & 2 & 0 & 1 & 1 & 2 & 2 & 3 & 1 \\
-2 & -1 & -1 & -1 & -2 & -1 & -1 & -1 & -2 & -2 & -1 & -1 & -1 & -2 & -2 & -2 & -1
\end{array}\right).
\end{align*}
}}
\subsection{Thinness of case C-56}
$B$ has infinite order, and the matrix is

\resizebox{\textwidth}{!}{
\parbox{\textwidth}{
\begin{align*}
&M=\left(\begin{array}{rrrrrrrrrrrrrrrrrrr}
0 & 0 & 0 & 0 & 0 & 1 & 1 & 1 & 1 & 2 & 2 & 3 & 3 & 4 & 4 & 5 & 5 \\
-2 & -1 & -1 & -1 & -1 & -2 & -1 & -1 & -1 & -5 & -3 & -7 & -5 & -9 & -7 & -11 & -10 \\
0 & 0 & 0 & 1 & 1 & 3 & 0 & 1 & 1 & 6 & 3 & 7 & 5 & 9 & 6 & 11 & 9 \\
-2 & -2 & -2 & -2 & -1 & -3 & -2 & -1 & -1 & -5 & -3 & -6 & -6 & -9 & -7 & -11 & -9 \\
0 & 0 & 0 & 0 & 1 & 3 & 0 & 1 & 1 & 5 & 2 & 7 & 5 & 10 & 7 & 11 & 9 \\
-1 & -2 & -1 & -2 & -1 & -2 & -2 & -2 & -1 & -3 & -1 & -4 & -2 & -5 & -3 & -5 & -4
\end{array}\right).
\end{align*}
}}
\subsection{Thinness of case C-57}
$B$ has infinite order, and the matrix is

\resizebox{\textwidth}{!}{
\parbox{\textwidth}{
\begin{align*}
 M= &\left(\begin{array}{rrrrrrrrrrrrrrrrrrrrrrrrrrrr}
0 & 0 & 0 & 0 & 0 & 0 & 1 & 1 & 1 & 2 & 2 & 3 & 3 & 4 & 4 & 5 & 5 & 5 & 5 & 7 & 21 & 28 \\
-18 & -15 & -4 & -4 & -4 & -1 & -16 & -2 & -1 & -5 & -3 & -7 & -5 & -9 & -7 & -20 & -13 & -11 & -10 & -8 & -63 & -63 \\
12 & 0 & 4 & 4 & 5 & 1 & 13 & 3 & 1 & 6 & 3 & 7 & 5 & 9 & 6 & 34 & 14 & 11 & 9 & -2 & 69 & 49 \\
0 & -15 & -3 & -1 & -3 & 0 & -6 & -3 & -1 & -5 & -3 & -6 & -6 & -9 & -7 & -59 & -14 & -11 & -9 & -5 & -84 & -43 \\
12 & 35 & 5 & 0 & 1 & 0 & 16 & 3 & 1 & 5 & 2 & 7 & 5 & 10 & 7 & 41 & 15 & 11 & 9 & 8 & 69 & 28 \\
-17 & -21 & -4 & -2 & -2 & -1 & -16 & -2 & -1 & -3 & -1 & -4 & -2 & -5 & -3 & -22 & -8 & -5 & -4 & -8 & -28 & -15
\end{array}\right).
\end{align*}
}}

\subsection{Thinness of case C-58}
$B$ has finite order, and the matrix is
$$M=\left(\begin{array}{rrrrrrrrrrrrrr}
0 & 0 & 0 & 1 & 1 & 1 & 2 & 2 & 3 & 3 & 4 & 4 & 5 & 5 \\
-1 & -1 & -1 & -2 & -2 & -1 & -5 & -3 & -7 & -5 & -9 & -7 & -11 & -10 \\
0 & 1 & 1 & 2 & 3 & 1 & 6 & 3 & 7 & 5 & 9 & 6 & 11 & 9 \\
-1 & -2 & -1 & -2 & -3 & -1 & -5 & -3 & -6 & -6 & -9 & -7 & -11 & -9 \\
0 & 1 & 1 & 1 & 3 & 1 & 5 & 2 & 7 & 5 & 10 & 7 & 11 & 9 \\
-1 & -1 & -1 & -1 & -2 & -1 & -3 & -1 & -4 & -2 & -5 & -3 & -5 & -4
\end{array}\right). $$

\section{Ping pong tables for \texorpdfstring{$\Sp(4)$}{Sp(4)} hypergeometric groups}\label{se:thinsp4}

In this section, we prove Theorem~\ref{thm:total-sp4} by establishing the thinness of cases 27, 35, 39 from Table~\ref{tab:sp4-all}. The solutions are represented by matrices, in the same way as in Section~\ref{se:thinlist}. The list of matrices is also accessible on GitHub~\cite{BDN21-s}.

\subsection{Thinness of case 27}
$B$ has infinite order, and the matrix is
$$ M=\left(\begin{array}{rrrrrr}
0 & 0 & 0 & 0 & 1 & 1 \\
-5 & -5 & -1 & -1 & -2 & -1 \\
0 & 4 & -1 & 0 & -3 & 1 \\
-1 & -5 & -1 & -5 & -2 & -1
\end{array}\right).$$

\subsection{Thinness of case 35}
$B$ has infinite order, and the matrix is
\begin{align*}
&M=\left(\begin{array}{rrrrrrrrrr}
0 & 0 & 0 & 1 & 1 & 1 & 2 & 2 & 3 & 3 \\
-1 & -1 & -1 & -2 & -2 & -1 & -5 & -3 & -7 & -6 \\
0 & 0 & 1 & 0 & 3 & 1 & 6 & 2 & 7 & 5 \\
-2 & -1 & -1 & -2 & -2 & -1 & -3 & -1 & -3 & -2
\end{array}\right).
\end{align*}

\subsection{Thinness of case 39}
$B$ has finite order, and the matrix is
$$ M=\left(\begin{array}{rrrrrrrrrr}
0 & 0 & 1 & 1 & 1 & 3 & 4 & 6 \\
-1 & -1 & -6 & -3 & -2 & -8 & -11 & -15 \\
0 & 1 & 4 & 3 & 2 & 9 & 11 & 16 \\
-1 & -1 & -3 & -3 & -1 & -6 & -6 & -9
\end{array}\right).$$

\section{Remarks on open cases}\label{se:open}
For the cases listed in Table~\ref{tab:open} we have not been able to show that $G=G_1\ast_{H} G_2$ (with the notation of Section~\ref{se:pingpong}). For some of these cases it might be that arithmeticity can be proven via the word method of \cite{bdss} or \cite{BDN22}, by searching even deeper than what has been done so far. In contrast, we checked that our ping pong approach must fail for all these cases and that this is not simply a matter of computation depth or of our choices for the rays that generate the cones.

Indeed, as mentioned in Section~\ref{se:construction}, the minimal choice for the starting cones $C_{0},D_{0}$ consists of a single ray $\mathbb{R}_{+}t_{0}$, where $t_{0}$ is a generator of the $1$-dimensional space $\mathrm{Im}(T-I)$. This is essentially a unique choice, up to swapping $C$ and $-C$. Then, we expand the cones through the procedure explained in Section~\ref{se:construction}: this is necessary in order to satisfy conditions \eqref{it:Bcheck} and \eqref{it:Tcheck}. In all cases of Table~\ref{tab:open}, after finitely many expansion steps we reach a pair of cones $C_{i},D_{i}$ for which condition~\eqref{it:disjoint} is violated. Therefore, in all these cases it is not possible to assume that the table half $X$ is given as the union of four cones $X^+=C\cup -C,\,X^-=EC\cup -EC$, with $E$ as in Section~\ref{se:pingpong}, and that it behaves under ping pong as depicted in Figure~\ref{fig:ppt}.

Still, a more complicated ping pong might work on these cases, and all three outcomes -- arithmetic, free, thin but not free -- are possible. One might conjecture that, like the resolved cases, all remaining cases are either arithmetic or (amalgamated) free products.

\section*{Acknowledgements}

The authors would like to thank for their hospitality the Technische Universit\"at Dresden, the Max Planck Institute f\"ur Mathematics (MPIM) in Bonn, and the Hebrew University of Jerusalem, where much of the work on this article was accomplished.

JB is supported by ERC Consolidator grant 681207 (codename GrDyAP, with A.~Thom as PI) and through the fellowship from MPIM, Bonn. DD is supported by the Israel Science Foundation Grants No. 686/17 and 700/21 of A.~Shalev, and the Emily Erskine Endowment Fund; he has been a postdoc at the Hebrew University of Jerusalem under A.~Shalev in 2020/21 and 2021/22. MN was supported by ERC Consolidator grant 681207 and DFG grant 281869850 (RTG 2229, ``Asymptotic Invariants and Limits of Groups and Spaces'').

In addition, the authors would like to thank the ``Research in Pairs'' program of Mathematisches Forschungsinstitut Oberwolfach (MFO) during which a significant progress in the  article was made.
\nocite{}
\bibliographystyle{abbrv}
\bibliography{BDN}

\begin{thebibliography}{10}

\bibitem{BDN21-s}
J.~Bajpai, D.~Dona, and M.~Nitsche.
\newblock Sage-{BDN}21.
\newblock SageMath code for ``Thin Monodromy in $\mathrm{Sp}(4)$ and
  $\mathrm{Sp}(6)$'' (Bajpai, Dona, Nitsche).
  \url{https://github.com/daniele-dona/Sage-BDN21}, 2021.

\bibitem{BDN22}
J.~Bajpai, D.~Dona, and M.~Nitsche.
\newblock {Arithmetic monodromy in $\mathrm{Sp}(2n)$}.
\newblock \texttt{arXiv:2209.07402}, 2022.

\bibitem{bdss}
J.~Bajpai, D.~Dona, S.~Singh, and S.~V. Singh.
\newblock Symplectic hypergeometric groups of degree six.
\newblock {\em J. Algebra}, 575:256--273, 2021.

\bibitem{BSS}
J.~Bajpai, S.~Singh, and S.~V. Singh.
\newblock Arithmeticity of some hypergeometric groups.
\newblock \texttt{arXiv:2201.08586}, 2022.

\bibitem{bhv}
B.~{Bekka}, P.~{de la Harpe}, and A.~{Valette}.
\newblock {\em {Kazhdan's {P}roperty ({T})}}, volume~11.
\newblock Cambridge: Cambridge University Press, 2008.

\bibitem{BH}
F.~Beukers and G.~Heckman.
\newblock Monodromy for the hypergeometric function {$_nF_{n-1}$}.
\newblock {\em Invent. Math.}, 95(2):325--354, 1989.

\bibitem{BT}
C.~Brav and H.~Thomas.
\newblock Thin monodromy in $\mathrm{Sp}(4)$.
\newblock {\em Compos. Math.}, 150(3):333--343, 2014.

\bibitem{CCJJV01}
P.~Cherix, M.~Cowling, P.~Jolissaint, P.~Julg, and A.~Valette.
\newblock {\em Groups with the {H}aagerup property: {G}romov's
  a-{T}-menability}, volume 197 of {\em Progress in mathematics}.
\newblock Springer, Basel (Switzerland), 2001.

\bibitem{DFH}
A.~S. Detinko, D.~L. Flannery, and A.~Hulpke.
\newblock Experimenting with symplectic hypergeometric monodromy groups.
\newblock {\em Exp. Math.}, to appear
  (https://doi.org/10.1080/10586458.2020.1780516).

\bibitem{DM06}
C.~F. Doran and J.~W. Morgan.
\newblock Mirror symmetry and integral variations of {H}odge structure
  underlying one parameter families of {C}alabi-{Y}au threefolds.
\newblock In N.~Yui, S.-T. Yau, and J.~D. Lewis, editors, {\em Mirror symmetry
  {V}, Proceedings of the {BIRS} Workshop on {C}alabi-{Y}au Varieties and
  Mirror Symmetry}, volume~38 of {\em AMS/IP Studies in Advanced Mathematics},
  pages 517--537. American Mathematical Society, International Press, 2006.

\bibitem{GMP}
B.~R. Greene, D.~R. Morrison, and M.~R. Plesser.
\newblock Mirror manifolds in higher dimension.
\newblock {\em Comm. Math. Phys.}, 173(3):559--597, 1995.

\bibitem{Jolissaint}
P.~Jolissaint.
\newblock Borel cocycles, approximation properties and relative property {T}.
\newblock {\em Ergodic Theory Dynam. Systems}, 20(2):483--499, 2000.

\bibitem{LTY}
B.~H. Lian, A.~Todorov, and S.-T. Yau.
\newblock Maximal unipotent monodromy for complete intersection {CY} manifolds.
\newblock {\em Amer. J. Math.}, 127(1):1--50, 2005.

\bibitem{LS77}
R.~C. Lyndon and P.~E. Schupp.
\newblock {\em Combinatorial group theory}.
\newblock Classics in Mathematics. Springer-Verlag, Berlin, 2001.
\newblock Reprint of the 1977 edition.

\bibitem{Sa14}
P.~Sarnak.
\newblock Notes on thin matrix groups.
\newblock In {\em Thin groups and superstrong approximation}, volume~61 of {\em
  Math. Sci. Res. Inst. Publ.}, pages 343--362. Cambridge Univ. Press,
  Cambridge, 2014.

\bibitem{S15S}
S.~Singh.
\newblock Arithmeticity of four hypergeometric monodromy groups associated to
  {C}alabi-{Y}au threefolds.
\newblock {\em Int. Math. Res. Not. IMRN}, 2015(18):8874--8889, 2015.

\bibitem{S17}
S.~Singh.
\newblock Arithmeticity of some hypergeometric monodromy groups in
  $\mathrm{Sp}(4)$.
\newblock {\em J. Algebra}, 473:142--165, 2017.

\bibitem{SS}
S.~Singh and S.~V. Singh.
\newblock Thinness of some hypergeometric groups in $\mathrm{Sp}(6)$.
\newblock \texttt{arXiv:2206.14159}, 2022.

\bibitem{SV}
S.~Singh and T.~N. Venkataramana.
\newblock Arithmeticity of certain symplectic hypergeometric groups.
\newblock {\em Duke Math. J.}, 163(3):591--617, 2014.

\end{thebibliography}

\end{document}